\documentclass[11pt]{amsart}
\usepackage{amsbsy}
\usepackage{amsmath}
\usepackage{amsfonts}
\usepackage{graphicx,epsfig,subfigure,psfrag}
\usepackage{color}

\usepackage{changes}

\DeclareMathOperator{\diver}{div}
\begin{document}

\newtheorem{theorem}{Theorem}[section]
\newtheorem{proposition}[theorem]{Proposition}
\newtheorem{lemma}[theorem]{Lemma}
\newtheorem{corollary}[theorem]{Corollary}
\newtheorem{definition}[theorem]{Definition}
\newtheorem{remark}[theorem]{Remark}
\newcommand{\tex}{\textstyle}
\numberwithin{equation}{section}
\newcommand{\ren}{\mathbb{R}^N}
\newcommand{\re}{\mathbb{R}}
\newcommand{\n}{\nabla}
\newcommand{\p}{\partial}
\newcommand{\iy}{\infty}
\newcommand{\pa}{\partial}
\newcommand{\fp}{\noindent}
\newcommand{\ms}{\medskip\vskip-.1cm}
\newcommand{\mpb}{\medskip}
\newcommand{\AAA}{{\bf A}}
\newcommand{\BB}{{\bf B}}
\newcommand{\CC}{{\bf C}}
\newcommand{\DD}{{\bf D}}
\newcommand{\EE}{{\bf E}}
\newcommand{\FF}{{\bf F}}
\newcommand{\GG}{{\bf G}}
\newcommand{\oo}{{\mathbf \omega}}
\newcommand{\Am}{{\bf A}_{2m}}
\newcommand{\CCC}{{\mathbf  C}}
\newcommand{\II}{{\mathrm{Im}}\,}
\newcommand{\RR}{{\mathrm{Re}}\,}
\newcommand{\eee}{{\mathrm  e}}
\newcommand{\LL}{L^2_\rho(\ren)}
\newcommand{\LLL}{L^2_{\rho^*}(\ren)}
\renewcommand{\a}{\alpha}
\renewcommand{\b}{\beta}
\newcommand{\g}{\gamma}
\newcommand{\G}{\Gamma}
\renewcommand{\d}{\delta}
\newcommand{\D}{\Delta}
\newcommand{\e}{\varepsilon}
\newcommand{\var}{\varphi}
\newcommand{\lll}{\l}
\renewcommand{\l}{\lambda}
\renewcommand{\o}{\omega}
\renewcommand{\O}{\Omega}
\newcommand{\s}{\sigma}
\renewcommand{\t}{\tau}
\renewcommand{\th}{\theta}
\newcommand{\z}{\zeta}
\newcommand{\wx}{\widetilde x}
\newcommand{\wt}{\widetilde t}
\newcommand{\noi}{\noindent}
\newcommand{\uu}{{\bf u}}
\newcommand{\xx}{{\bf x}}
\newcommand{\yy}{{\bf y}}
\newcommand{\zz}{{\bf z}}
\newcommand{\aaa}{{\bf a}}
\newcommand{\cc}{{\bf c}}
\newcommand{\jj}{{\bf j}}
\newcommand{\ggg}{{\bf g}}
\newcommand{\UU}{{\bf U}}
\newcommand{\YY}{{\bf Y}}
\newcommand{\HH}{{\bf H}}
\newcommand{\GGG}{{\bf G}}
\newcommand{\VV}{{\bf V}}
\newcommand{\ww}{{\bf w}}
\newcommand{\vv}{{\bf v}}
\newcommand{\hh}{{\bf h}}
\newcommand{\di}{{\rm div}\,}
\newcommand{\ii}{{\rm i}\,}
\def\I{{\mathbf{I}}}
\newcommand{\inA}{\quad \mbox{in} \quad \ren \times \re_+}
\newcommand{\inB}{\quad \mbox{in} \quad}
\newcommand{\inC}{\quad \mbox{in} \quad \re \times \re_+}
\newcommand{\inD}{\quad \mbox{in} \quad \re}
\newcommand{\forA}{\quad \mbox{for} \quad}
\newcommand{\whereA}{,\quad \mbox{where} \quad}
\newcommand{\asA}{\quad \mbox{as} \quad}
\newcommand{\andA}{\quad \mbox{and} \quad}
\newcommand{\withA}{,\quad \mbox{with} \quad}
\newcommand{\orA}{,\quad \mbox{or} \quad}
\newcommand{\atA}{\quad \mbox{at} \quad}
\newcommand{\onA}{\quad \mbox{on} \quad}
\newcommand{\ef}{\eqref}
\newcommand{\mc}{\mathcal}
\newcommand{\mf}{\mathfrak}

\newcommand{\ssk}{\smallskip}
\newcommand{\LongA}{\quad \Longrightarrow \quad}
\def\com#1{\fbox{\parbox{6in}{\texttt{#1}}}}
\def\N{{\mathbb N}}
\def\A{{\cal A}}
\newcommand{\de}{\,d}
\newcommand{\eps}{\varepsilon}
\newcommand{\be}{\begin{equation}}
\newcommand{\ee}{\end{equation}}
\newcommand{\spt}{{\mbox spt}}
\newcommand{\ind}{{\mbox ind}}
\newcommand{\supp}{{\mbox supp}}
\newcommand{\dip}{\displaystyle}
\newcommand{\prt}{\partial}
\renewcommand{\theequation}{\thesection.\arabic{equation}}
\renewcommand{\baselinestretch}{1.1}
\newcommand{\Dm}{(-\D)^m}


 \title[Logistic equation]{\bf Interface logistic problems: large diffusion and singular perturbation results}

\author{Pablo \'Alvarez-Caudevilla, Cristina Br\"{a}ndle, M\'onica Molina-Becerra and Antonio Su\'arez}

\address{Universidad Carlos III de Madrid, Spain} \email{pacaudev@math.uc3m.es, cbrandle@math.uc3m.es}
\address{Universidad de Sevilla, Spain} \email{monica@us.es, suarez@us.es}



\date{\today}

\parskip5pt

\begin{abstract}
In this work we consider an interface logistic problem where two populations live in two different regions, separated by a membrane or interface where it happens  an
interchange of flux. Thus, the two populations only interact or are coupled through such a membrane where we impose the so-called Kedem-Katchalsky boundary conditions.
For this particular scenario we analyze the existence and uniqueness of positive solutions
depending on the parameters involve in the system, obtaining interesting results where one can see for the first time the effect of the membrane under such boundary conditions. To do so, we first ascertain the
asymptotic behaviour of several linear and nonlinear problems for which we include a diffusion coefficient and analyse the behaviour of the solutions when such a diffusion parameter goes to zero or infinity. Despite their
own interest, since these asymptotic results have never been studied before, they will be crucial in analyzing the existence and uniqueness for the main interface logistic problems under analysis.
Finally, we apply such an asymptotic
analysis to characterize the existence of solutions in terms of the growth rate of the populations, when both populations possess the same growth rate and, also, when they depend on different parameters.

\end{abstract}

\maketitle

\section{Introduction}
\begin{center}
\begin{figure}[ht!]
\centering{
\includegraphics[scale=0.35]{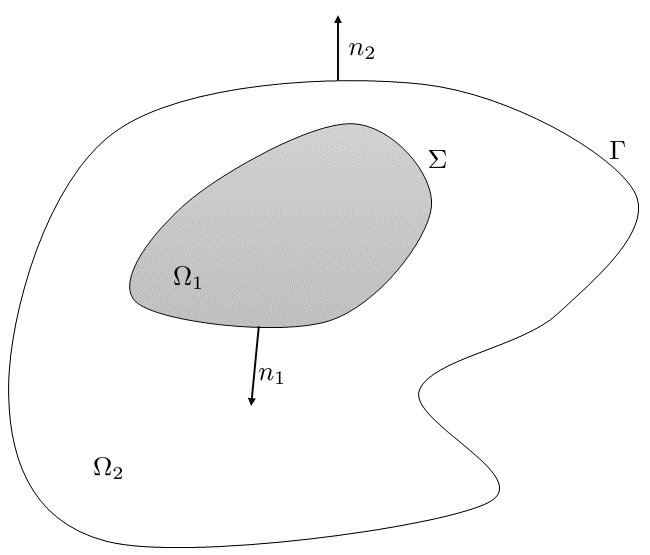}
\caption{A possible configuration of the domain $\Omega=\Omega_1\cup \Omega_2 \cup \Sigma$.}}
\label{figure}
\end{figure}
\end{center}

Let $\Omega$ be a bounded domain of $\mathbb{R}^N$ split by
\begin{equation}
\label{eq:domain}
\Omega=\Omega_1\cup \Omega_2 \cup \Sigma,
\end{equation}
such that $\Omega_i$, $i=1,2$ are two subdomains, and
$\Sigma=\partial \Omega_1$ stands for an internal boundary, while $\Gamma=\partial\Omega_2\setminus\Sigma$ represents the exterior boundary; see
Figure \ref{figure} where we have illustrated an example of $\Omega$ showing such a configuration.
 The main goal of this paper is to study existence and uniqueness of positive solutions as well as the limit behaviour when the diffusion coefficient tends to zero or infinity  of different logistic systems of the form
\begin{equation}
  \label{eq:general-system}
  \left\{\begin{array}{ll}
-d\Delta u_{1} =u_{1}(\beta_1(x)-\alpha_1 (x)u_{1})\quad   \hbox{in } \Omega_1,\\
-d\Delta u_{2} =u_{2}(\beta_2(x)-\alpha_2 (x)u_{2})\quad   \hbox{in } \Omega_2,
 \end{array}\right.
\end{equation}
where the equations are coupled only through the so-called interface boundary conditions
\begin{equation}
  \label{eq:KK-condition}
  \left\{
  \begin{array}{l}
  \partial_{\bf n_1}u_1=\gamma_1(u_2-u_1)\\
  \partial_{\bf n_1}u_2=\gamma_2(u_2-u_1)
  \end{array}
  \right. \quad \hbox{on } \Sigma.
\end{equation}
Here ${\bf n_i}$ is the outward normal vector to $\Omega_i$, and ${\bf n_1}=-{\bf n_2}$ on $\Sigma$ so that
$$
\partial_{\bf n_1} u_2=\gamma_2(u_2-u_1) \Longleftrightarrow \partial_{\bf n_2}u_2=\gamma_2(u_1-u_2).
$$
We also impose homogeneous Neumann conditions on the outside boundary.
\begin{equation}
  \label{eq:neumann-condition}
  \partial_{\bf n_2}u_2=0 \quad \hbox{on}\quad \Gamma.
\end{equation}
The parameter $d>0$ stands for the diffusion coefficient and $\gamma_i>0$ represents the membrane
permeability effect of the barrier or interface $\Sigma$ for
the movement of flux from one subdomain to the other. The coefficients $\alpha_i(x)$
 are two positive  continuous functions in  ${C}(\overline{\Omega_i})$ that, from a Biological point of view, measure the crowding effects of
the populations in the corresponding subdomains of $\Omega_i$ and $\beta_i$, which are also continuous,  may change sign, but $\max_{x\in\Omega_i}\beta_i(x)>0$.
 Moreover, we will also be interested in studying problems in which $\beta_i(x)$ is constant, so that $\beta_i(x)=\lambda_i>0$ stand for the intrinsic growth rate of the populations. In these cases we will take the diffusion to be $1$ and the interest will be focused on moving the parameters $\lambda_i$. Although the two limit behaviors that we consider here ($d$ to 0 or infinity and $\lambda_i$ to infinity) may seem a priori very different, they are, as we shall see, strongly related and to analyze the second one it is essential to know the first one.

This type of interface logistic systems has been of increasing interest in recent years, since they are used to model different problems in physical or biological contexts. In particular, the case in which one assumes  constant growth rates
($\beta_i(x)$ constant and the same on both sides of the domain $\Omega$, $\beta_1=\beta_2=\lambda$) has been deeply studied during the last years, having nowadays
a wide knowledge about the behaviour of its solutions.  One can see such an analysis of the behaviour of its solutions in
the works \cite{PAC-Bran, Chaplain, CiavolellaPaerthame, Ciavolella, Suarezetal}, where many applications
and interest in other sciences are described, despite its own mathematical interest.

The boundary conditions \eqref{eq:KK-condition} describe the flow through the common boundary $\Sigma$ for the subdomains. Those boundary conditions
are compatible with mass conservation and energy dissipation, and are called Kedem-Katchalsky membrane conditions.
These Kedem-Katchalsky boundary conditions were introduced in a thermodynamic context in \cite{Ked-Kat}. However, their application has been also studied for Biological problems
later on by Quarteroni et al \cite{Quarteroni} in the analysis of the dynamics of the solute in a vessel and an arterial wall.
The conservation of mass provides us with
the flux continuity, and the dissipation principle implies that the $L^2$-norm
of the solution is decreasing in time. This last property is obtained assuming that the flux as proportional to the difference of densities on the interface region or barrier; see \cite{Ciavolella} for a discussion of how to deduce such boundary conditions under those physical assumptions.  Moreover, note that here we assume different membrane permeability constants of proportionality, $\gamma_i$, depending on the direction of the flow.

 From a mathematical point of view we would like to mention the works of Chen \cite{Chen} where he studied the effects of the barrier on the global dynamics of a parabolic problem with similar boundary conditions to \eqref{eq:KK-condition}, establishing
the existence of weak solutions, a comparison principle and some stability results. Also, for the parabolic problem, in very recent works Ciavolella and Perthame in \cite{CiavolellaPaerthame} obtained several regularity results adapting the classical
$L^1$ theory for parabolic reaction-diffusion systems. With respect to the elliptic models assuming such boundary conditions, up to our knowledge, there are only a couple of works \cite{PAC-Bran, Suarezetal, WangSu} where these elliptic models were under analysis
assuming some spatial heterogeneities with symmetric boundary conditions and  the existence and uniqueness of positive solutions with non-symmetric boundary conditions as \eqref{eq:KK-condition}. Finally, we must point out that the results of this paper work for more general geometric configurations as those shown by Pfl\"{u}ger \cite{Pfluger1996}.

\noindent{\bf Notation.} First, we fix some notations that will be used throughout the paper. For convenience, we will use the subscript $1$, respectively $2$, to denote objects (functions, parameters...) defined in $\Omega_1$ (resp. $\Omega_2$) and each time we write the subscript $i$ we mean $i=1,2$ and we will not mention it anymore unless it is necessary.
We write bold letters to denote vectors ${\bf w}=(w_1,w_2)^T$ with $w_i$ defined in $\Omega_i$. Thus, the first entry is defined in $\Omega_1$, or defined equal to zero in $\Omega_2$, and the second one in $\Omega_2$.
For example, for the solution ${\bf u}$ of the system~\eqref{eq:general-system} one has that $u_i:\Omega_i\to\mathbb{R}$ and ${\bf u}(x)$ will stand for $(u_1(x), 0)^T$ if $x\in\Omega_1$ and $(0,u_2(x))^T$ if $x\in\Omega_2$. By abuse of notation we will use  $\nabla {\bf u}$ to denote  $(\nabla u_1,\nabla u_2)^T$. Also, we will write ${\bf w}>0$ if $w_1>0 $ in $\Omega_1$ and $w_2>0$ in $\Omega_2$. 
Furthermore, in order to simplify the notation we write the boundary conditions \eqref{eq:KK-condition}--\eqref{eq:neumann-condition} as
\begin{equation}
\label{eq:boundary_mem}
\mathcal{I}({\bf u})=0 \quad\mbox{on $\Sigma\cup \Gamma$.}
\end{equation}
Let $\mathcal{U}$ be regular and bounded domain $\mathcal{U}\subset \mathbb{R}^N$. Given $f\in C(\overline{\mathcal{U}})$ we denote
$$
f_L:=\min_{x\in \overline{\mathcal{U}}}f(x),\qquad f_M:=\max_{x\in \overline{\mathcal{U}}}f(x)\qquad \text{and} \qquad
f_+=\max\{0,f\}.
$$

\noindent{\bf Main results.}
 As said, the main goal of the paper is to study the existence and uniqueness of positive solutions of problem~\eqref{eq:general-system}--\eqref{eq:neumann-condition} under different scenarios, namely when $\beta_i(x)$ are two different
 functions that are allowed to change sign and when they are different constants $\lambda_i$. We are also concerned about the limit when the diffusion coefficient $d$ degenerates or, on the contrary, becomes very large.
To this aim, we must state before some important auxiliary results concerning the asymptotic behaviour (when $d\to 0, \infty)$ of linear and non-linear logistic scalar problems. This particular analysis on the effect of large diffusion will be crucial in the sequel to determine the existence of positive solutions of problem ~\eqref{eq:general-system}--\eqref{eq:neumann-condition}. Nevertheless, they have their own interest and
mathematical relevance since they have never been studied before.

In particular, in Section~\ref{Sec:scalar} we consider the logistic scalar version of~\eqref{eq:general-system}
$$
-d\Delta u =u(\beta(x)-\alpha(x)u) \quad \hbox{in } D,
$$  where $D\subset \mathbb{R}^N$ is a regular, bounded domain such that $
\partial D=\Gamma_1\cup \Gamma_2.
$
On each part of the boundary, we will consider the boundary conditions
$$
{\mathcal{B}}u:=
{\partial_{\bf n}} u +g_i u \quad  \mbox{on $\Gamma_i$,}
$$
first of homogeneous and then of non-homogeneous type, where $g_i$ are regular functions and ${\bf n}$ is a nowhere tangent vector field .
We study the behaviour of the unique positive solution  of these problems  when the diffusion parameter tends to zero and to infinity. We will see that, when $d\to 0$
the solution to the limit problem exists and its behaviour will depend on the relation between $\alpha$ and $\beta$. On the other hand, when $d\to\infty$ the existence and
behaviour of the limit solution depends on the sign of $\sigma_1^D[-\Delta;{\mathcal{B}}]$, the principal eigenvalue of the problem $-\Delta u=\sigma u$ under the boundary conditions given by $\mathcal{B}$.
Related to that results, and as a previous step, we also analyse  some monotonicity properties and convergence of the principal eigenvalue of the associated linear eigenvalue problem with equation
$$
-d\Delta \varphi  -\beta(x) \varphi=\mu \varphi\quad \hbox{in } D.
$$

Subsequently, in Section~\ref{Sec:Interface_eigenvalue} we ascertain some similar asymptotic results to those in Section~\ref{Sec:scalar}, now for an interface eigenvalue system of the form
$$
-d\Delta \varphi_i + c_i(x)
\varphi_i=\nu \varphi_i\quad  \hbox{in } \Omega_i,\\
$$
with the interface boundary conditions given by~\eqref{eq:boundary_mem}.
We obtain the convergence of the principal
eigenvalue associated with such an interface eigenvalue problem when the parameter $d$ goes to zero and also to infinity. It is worth mentioning that, contrary to what happens for the scalar eigenvalue problem studied in Section~\ref{Sec:scalar}, the only principal eigenvalue of $-\Delta u_i=\sigma u_i$,  in $\Omega_i$,  under the membrane boundary conditions~\eqref{eq:boundary_mem}, is $\Lambda_1(-\Delta,-\Delta)=0$. The associate principal eigenfunctions are constant. This fact impacts on  the limit behaviour, when $d\to\infty$ of the principal eigenvalue to the interface eigenvalue problem, see Propositions~~\ref{prop.scalar.eigenvalue} and \ref{prop.eigenvalue.d.to.infinity}; showing this interesting
effect of the membrane on such limit behaviours.

The results obtained in Sections~\ref{Sec:scalar} and~\ref{Sec:Interface_eigenvalue} will be crucial in the sequel analysis of the nonlinear interface logistic system \eqref{eq:general-system}--\eqref{eq:neumann-condition}. In particular, in Section~\ref{Sec:interface_logistic},
 using some upper bounds for the solutions, already proved in \cite{Suarezetal}, we find the convergence of the solutions of \eqref{eq:general-system}--\eqref{eq:neumann-condition} when either $d\to 0$ or $d\to \infty$.  Since these limit behaviours depend strongly on the eigenvalue problem studied in Section~\ref{Sec:Interface_eigenvalue}, we observe again that, in case $d\to 0$  the limit solutions in the scalar case and in the interface problem are equivalent. On the contrary, when $d\to\infty$ the behavior differs radically.

Finally, in Section~\ref{Sec:interface_growth_rates}
we focus on the interface problem \eqref{eq:general-system}--\eqref{eq:neumann-condition} when $d=1$ and $\beta_i(x)$ are taken as real parameters $\lambda_i$, both in the case in which $\lambda_1=\lambda_2$
(similar to the previous existing works) and when the parameters differ (not studied before).

As mentioned above this interface logistic problem with constant growth rates  is nowadays considered the more classical logistic interface problem, thoroughly study during the last years, for example in
the works \cite{PAC-Bran, Ciavolella, Suarezetal}. However, the case when the growth rates are different, especially with two different parameters, has not been analysed before.
In particular, we apply the asymptotic results obtained in the previous sections, while moving the diffusion parameter $d$, to ascertain the existence and uniqueness of solutions.
 More precisely,  we characterize the existence and uniqueness of solutions depending on the size of $\lambda_i$.

First, when $\lambda_1=\lambda_2=\lambda$, as direct consequence of the asymptotic behaviour of the solutions when the parameter $d$ goes to zero or infinity performed in Section~\ref{Sec:interface_logistic} we arrive at a characterization of the solutions {$\bf  \Theta$} of the problem. Namely, taking $\lambda$ to be $1/d$, we find a for relation ${\theta_{i}}/{\lambda}$, both when $\lambda\to 0$ and $\lambda\to\infty$ in terms of $\alpha_i$ and the size of the domains $\Omega_i$.  On the other hand, when $\lambda_1\neq \lambda_2$ the solution {$\bf \Theta_\lambda$} of the problem exists and is unique if  either $\lambda_2\geq \sigma_2$ and $\lambda_1\in \mathbb{R}$, or
$\lambda_2<\sigma_2$ and $\lambda_1>\mathcal{H}(\lambda_2)$, where $\sigma_1,\sigma_2>0$ and the regular and decreasing function $\mathcal{H}$ are established in Section~\ref{Sec:interface_growth_rates} below. So that for a fixed value of $\lambda_2$ we study the behaviour for the solution when $\lambda_1\to\pm\infty$.  When $\lambda_1\to\infty$, the \lq\lq inside\rq\rq\ component of the solution, $\theta_1$, blows-up in the whole domain $\Omega_1$, whereas  $\theta_2$ blows up on the inner boundary $\Sigma$ and remains bounded in $\Omega_2$ converging to the so-called large solution. On the other hand, when $\lambda_1\to-\infty$, $\theta_1$ vanishes and $\theta_2$ tends to the unique positive solution of a classical logistic equation with homogeneous boundary data.

\section{Preliminaries - Scalar problems}
\label{Sec:scalar}

In this section we consider two scalar problems, namely an eigenvalue problem and a general logistic equation and show several auxiliary results that will be used in the sequel for proving the main results obtained
in the paper.  Throughout the section we will consider a domain $D\subset \mathbb{R}^N$ which is regular and bounded, such that
$$
\partial D=\Gamma_1\cup \Gamma_2,
$$
where $\Gamma_1$ and $\Gamma_2$ are two disjoint open and closed subsets of the boundary, consistently, in some way, with the configuration shown in Figure \ref{figure}.
Note that, it is possible that some $\Gamma_i=\emptyset$ for some $i=1,2$. On the boundary we will impose a Robin-type boundary condition of the form
\begin{equation}
\label{eq:BC_eign_prob}
{\mathcal{B}}\varphi:=
{\partial_{\bf n}} \varphi +g_i \varphi \quad  \mbox{on $\Gamma_i$,}
\end{equation}
where ${\bf n}$ is a nowhere tangent vector field and $g_i\in C(\Gamma_i)$.

\subsection{Scalar eigenvalue problem}
Let us consider first the eigenvalue problem
\begin{equation}
\label{eq:eign_prob_logistic}
 \left\{\begin{array} {l@{\quad}l} -d\Delta \varphi + c(x)
\varphi=\mu \varphi& \hbox{in } D,\\ \mathcal{B} \varphi={\bf 0}& \hbox{on } \partial D.
\end{array}\right.
\end{equation}
 Here  $c\in \mathcal{C}(D)$. We denote the principal eigenvalue of \eqref{eq:eign_prob_logistic}  as
\begin{equation}
\label{eq:princ_eign}
\sigma_1^D[-d\Delta+c;{\mathcal{B}}]=\sigma_1^D[-d\Delta+c;\mathcal{N}+g_1;\mathcal{N}+g_2]
\end{equation}
with $\mathcal{N}$ standing for the normal derivative shown in the general boundary conditions \eqref{eq:BC_eign_prob}.

First we also include a continuous dependence of the principal eigenvalue with respect
to different coefficients of (\ref{eq:eign_prob_logistic}).
\begin{proposition}
Let $\mathcal{R}:=L^\infty(D)\times C(\Gamma_1)\times C(\Gamma_2)$. Then,
\begin{enumerate}
\label{proper}
\item The map $(d,c,g_1,g_2)\in (0,+\infty)\times \mathcal{R}  \mapsto    \sigma_1^D[-d\Delta+c;{\mathcal{B}}]\in \mathbb{R}$
is continuous.
\item The map $(c,g_1,g_2)\in \mathcal{R} \mapsto  \sigma_1^D[-d\Delta+c;{\mathcal{B}}]\in  \mathbb{R}$
is increasing.
\end{enumerate}
\end{proposition}

\begin{proof}
The proof can be done following the arguments shown in \cite{Santi}.
\end{proof}

Now we deal with the limit of the principal eigenvalue when the diffusion coefficient $d$ tends to zero or infinity. The next result ascertains the behaviour as $d$ goes to zero and it follows from Theorem 4.1 in \cite{sergiojlg} (see also Proposition 1.3.6 in \cite{lam-lou}).

\begin{proposition}
\label{prop.scalar.limit.0}
Let $\sigma_1^D[-d\Delta+c;{\mathcal{B}}]$ be the principal eigenvalue of~{\rm \eqref{eq:eign_prob_logistic}}. Then:
$$
\lim_{d\to 0}\sigma_1^D[-d\Delta+c;{\mathcal{B}}]=c_L.
$$
\end{proposition}

The limit as $d$ tends to infinity is analysed in the next result.
\begin{proposition}
\label{prop.scalar.eigenvalue}
Let $\varphi_*$ be the principal normalized eigenfunction, $\|\varphi_*\|_\infty=1$, associated with $\sigma_1^D[-\Delta;{\mathcal{B}}]$.
Then
\begin{equation}
\label{eq:lim_princ_eign}
\lim_{d\to \infty }\sigma_1^D[-d\Delta+c;{\mathcal{B}}]=
\left\{
\begin{array}{ll}
+\infty & \mbox{\quad if $\sigma_1^D[-\Delta;{\mathcal{B}}]>0$,}\\
-\infty & \mbox{\quad  if $\sigma_1^D[-\Delta;{\mathcal{B}}]<0$,}\\
\displaystyle\left(\int_D c(x)\varphi_*^2 \right)\left(\displaystyle\int_D\varphi_*^2\right)^{-1} & \mbox{\quad  if $\sigma_1^D[-\Delta;{\mathcal{B}}]=0$.}\\
\end{array}
\right.
\end{equation}
\end{proposition}

\begin{proof}
We observe that for the minimum value of the potential function $c(x)$ it follows that
 $$
 \sigma_1^D[-d\Delta+c;{\mathcal{B}}]\geq\sigma_1^D[-d\Delta+c_L;{\mathcal{B}}]=\sigma_1^D[-d\Delta;{\mathcal{B}}]+c_L=d\sigma_1^D[-\Delta;{\mathcal{B}}]+c_L.
 $$
Analogously, for the maximum value it yields
$$
 \sigma_1^D[-d\Delta+c;{\mathcal{B}}]\leq\sigma_1^D[-d\Delta;{\mathcal{B}}]+c_M=d\sigma_1^D[-\Delta;{\mathcal{B}}]+c_M.
$$
Hence, combining both of them
\begin{equation}
\label{alavez}
 d\sigma_1^D[-\Delta;{\mathcal{B}}]+c_L\leq  \sigma_1^D[-d\Delta+c;{\mathcal{B}}]\leq d\sigma_1^D[-\Delta;{\mathcal{B}}]+c_M.
 \end{equation}
If  $\sigma_1^D[-\Delta;{\mathcal{B}}]\neq 0$, since it is also bounded, we deduce the first two parts of the limit \eqref{eq:lim_princ_eign}.

 Now, assume that $\sigma_1^D[-\Delta;{\mathcal{B}}]=0$. Then, due to (\ref{alavez})
$$
c_L\leq\sigma_1^D[-d\Delta+c;{\mathcal{B}}]\leq c_M,
$$
so that $\sigma_1^D[-d\Delta+c;{\mathcal{B}}]$ is bounded.

We now denote $\varphi_d$ as the positive eigenfunction associated with  $\sigma_1^D[-d\Delta+c;{\mathcal{B}}]$, normalized by $\|\varphi_d\|_\infty=1$, that is
\begin{equation}
\label{fid}
-d\Delta \varphi_d+c(x)\varphi_d=\sigma_1^D[-d\Delta+c;{\mathcal{B}}]\varphi_d,\quad \mbox{in $D$,}\quad {\mathcal{B}}\varphi_d={\bf 0} \quad \mbox{on $\partial D$.}
\end{equation}
We can choose $K>0$ such that $\sigma_1^D[-\Delta+K;{\mathcal{B}}]>0$, which leads to the problem
 $$
(-\Delta+K) \varphi_d=\left(\frac{\sigma_1^D[-d\Delta+c;{\mathcal{B}}]-c(x)}{d}+K\right)\varphi_d\quad \mbox{in $D$,}\quad {\mathcal{B}}\varphi_d={\bf 0} \quad \mbox{on $\partial D$.}
$$
Hence, by elliptic regularity $\|\varphi_d\|_{W^{2,p}}\leq C$ for some positive constant $C$ and any $p>1$. Therefore,
$$
\varphi_d\to\varphi_*\quad\mbox{in $C^1(\overline D)$ as $d\to \infty$, with $\varphi_*>0$, $\|\varphi_*\|_\infty=1$,}
$$
where $\varphi_*$ is a positive solution of
\begin{equation}
\label{fiestre}
 -\Delta \varphi_*=0\quad \mbox{in $D$,}\qquad {\mathcal{B}}\varphi_*={\bf 0} \quad \mbox{on $\partial D$.}
\end{equation}
 Finally, multiplying (\ref{fid}) by $\varphi_*$, and (\ref{fiestre}) by $d\varphi_d$, integrating by parts and subtracting both expression, we get to
$$
\int_D c(x)\varphi_d \varphi_*=\sigma_1^D[-d\Delta+c;{\mathcal{B}}]\int_D\varphi_d \varphi_*.
$$
Consequently, passing to the limit as $d\to \infty$ we find that
$$
\lim_{d\to \infty}\sigma_1^D[-d\Delta+c;{\mathcal{B}}]=\lim_{d\to\infty} \frac{\displaystyle\int_D c(x)\varphi_d \varphi_*}{\displaystyle\int_D\varphi_d \varphi_*}= \frac{\displaystyle\int_D c(x)\varphi_*^2}{\displaystyle\int_D\varphi_*^2}.
$$
\end{proof}
\begin{remark}
The particular case ${\mathcal{B}}={\mathcal{N}}$ has been studied in Proposition {\rm 1.3.19} of \cite{lam-lou}. In such case, $\sigma_1^D[-d\Delta+c;{\mathcal{N}}]=0$ and $\varphi_*=K$, $K$ constant. Hence,  $\lim_{d\to \infty}\sigma_1^D[-d\Delta+c;{\mathcal{N}}]=\frac{1}{|D|}\int_D c(x)$.
 \end{remark}
\subsection{Scalar logistic equation}
Now consider the scalar version of the general problem~\eqref{eq:general-system}, the classical logistic equation
\begin{equation}
\label{eq:classic_logistic}
 \left\{\begin{array} {l@{\quad}l} -d\Delta u =u(\beta(x)-\alpha(x)u)& \hbox{in } D,\\
 \mathcal{B}u={\bf 0}& \hbox{on } \partial D,
\end{array}\right.
\end{equation}
where, as mentioned in the Introduction,  $\beta\in C(\overline\Omega)$ might change sign, with $\beta_M>0$ and  $\alpha\in C(\overline\Omega)$ is such that $\alpha(x)>0$ for $x\in\overline\Omega$.

The next existence result for problem \eqref{eq:classic_logistic} follows from \cite{Santi2}
\begin{proposition}
\label{prop:existence!.logistic}
There exists a positive solution of {\rm(\ref{eq:classic_logistic})} if and only if
$$
\sigma_1^D[-d\Delta-\beta;{\mathcal{B}}]<0.
$$
When the positive solution exists, it is  unique.
\end{proposition}
Observe that thanks to Proposition~{\rm \ref{prop.scalar.limit.0}} we get
$$
\lim_{d\to 0}\sigma_1^D[-d\Delta-\beta;{\mathcal{B}}]=(-\beta)_L=-\beta_M<0,
$$
and then one can say that {\rm(\ref{eq:classic_logistic})} possesses a unique positive solution for $d$ small.

\begin{proposition}
\label{logisnormal}
  Let ${\tilde u}_d$ be the unique positive solution of {\rm(\ref{eq:classic_logistic})}. Then
  \begin{equation*}
\label{limidcero}
\lim_{d\to 0}{\tilde u}_d=\frac{\beta_+(x)}{\alpha(x)}\quad\mbox{uniformly in $\overline\Omega$.}
\end{equation*}
\end{proposition}

\begin{remark}
 The proof of this result follows from \cite{sergiojlg}, see also Theorem {\rm 5.2.5} in the particular case of Neumann boundary conditions.
\end{remark}

Our next aim is establishing  conditions  for the existence (or not) of solutions to~\eqref{prop.scalar.eigenvalue} in the limit case  $d\to\infty$. These existence results depend on the sign of $\sigma_1^D[-\Delta;{\mathcal{B}}]$. We will need to apply a version of the Picone's identity, that we include for the sake of completion.

\begin{lemma}
\label{lemma:Picone}
Let $q>N$ and $u,v\in W^{2,q}(D)$
be arbitrary functions in the domain $D$, such that $u$ and $v$ satisfy $\mathcal{B}u=\mathcal{B}v=0$  on
$\partial D$ and $\frac{v}{u}\in \mc{C}^1(D)\cap \mc{C}(\overline D)$.
Then, for any function
$f:\mathbb{R} \longrightarrow \mathbb{R}$ of class $\mc{C}^1$ the following Picone's identity holds:
\begin{equation*}
\label{216}
	\int_{D} f\left(\frac{v}{u}\right)
	(v \Delta u- u \Delta v)
	=\int_{D} f'\left(\frac{v}{u}\right) u^2 \left| \nabla \left(\frac{v}{u} \right) \right|^2.
\end{equation*}

\begin{proof}
Rearranging terms, it follows that
\begin{align*}
 f\left(\frac{v}{u}\right) (v \Delta u- u \Delta v) & =
  f\left(\frac{v}{u}\right) \left[v\diver(\nabla u)  - u\diver(\nabla v)\right] \\
  &=\diver\left(f\left(\frac{v}{u}\right) v \nabla u -f\left(\frac{v}{u}\right) u \nabla v \right)-\nabla\left(f\left(\frac{v}{u}\right) v\right)\nabla u
  +\nabla \left(f\left(\frac{v}{u}\right) u\right)\nabla v.
\end{align*}
Taking the first term on the right hand side, after integrating in $D$ we find that, applying the Divergence Theorem,
\begin{align*}
\int_D \diver\left(f\left(\frac{v}{u}\right) v \nabla u -f\left(\frac{v}{u}\right) u \nabla v \right)
= \int_{\partial D} f\left(\frac{v}{u}\right)  \left(v \frac{\partial u}{\partial {\bf n}} - u \frac{\partial v}{\partial {\bf n}}\right)=0,
\end{align*}
where the last equality is obtained after applying the boundary conditions.
On the other hand with respect to the second term differentiating and rearranging terms
we arrive at
$$
-\nabla \left(f\left(\frac{v}{u}\right) v\right)\nabla u
  + \nabla\left(f\left(\frac{v}{u}\right) u\right)\nabla v=
  \nabla\left(f\left(\frac{v}{u}\right) \right) \left( -v\nabla u+u\nabla v\right)=f'\left(\frac{v}{u}\right) \left|\nabla\left(\frac{v}{u}\right)\right|^2u^2,
  $$
completing the proof.
\end{proof}

\end{lemma}
\begin{proposition}
\label{dinfilogistica} Let ${\tilde u}_d$ be the unique positive solution of {\rm(\ref{eq:classic_logistic})}. Then
$$
\lim_{d\to\infty} \tilde{u}_d(x)= \left\{
\begin{array}{ll}
+\infty & \mbox{if $\sigma_1^D[-\Delta;{\mathcal{B}}]<0$,}\\
\tilde{u}_\infty & \mbox{if $\sigma_1^D[-\Delta;{\mathcal{B}}]=0$ and $\displaystyle  \int_D \beta(x)\varphi_*^2>0$,}
\end{array}
\right.
$$
where $\tilde{u}_\infty$ is a positive solution of the problem
\begin{equation}
\label{eq:inft_u}
-\Delta \tilde{u}_\infty=0\quad\mbox{in $D$,}\quad \mathcal{B}\tilde{u}_\infty={\bf 0}\quad\mbox{on $\partial D$.}
\end{equation}
\end{proposition}

\begin{proof}
  Assume $\sigma_1^D[-\Delta;{\mathcal{B}}]<0$. Take $\varphi_*$ the positive eigenfunction associated with $\sigma_1^D[-\Delta;{\mathcal{B}}]$ such that $\|\varphi_*\|_\infty=1$. Then, $\varepsilon \varphi_*$ is subsolution of   (\ref{eq:classic_logistic}) provided that
$$
\varepsilon \alpha(x)\varphi_*\leq \beta(x)-d\sigma_1^D[-\Delta;{\mathcal{B}}].
$$
Then, for $d$ large, we obtain
$$
\frac{\beta_L-d\sigma_1^D[-\Delta;{\mathcal{B}}]}{\alpha_M}\varphi_*(x)\leq \tilde{u}_d(x).
$$
Finally, since $\sigma_1^D[-\Delta;{\mathcal{B}}]<0$ and $\varphi_*(x)>0$ for all $x\in\overline\Omega$ we conclude that $\displaystyle\lim_{d\to\infty} \tilde{u}_d(x)=+\infty$.

On the other hand, if $\sigma_1^D[-\Delta;{\mathcal{B}}]=0$, due to Proposition  \ref{prop.scalar.eigenvalue} we have that
\begin{equation}
  \label{eq:noexistence.sigma0}
\lim_{d\to\infty }\sigma_1^D[-d\Delta-\beta;{\mathcal{B}}]=-\frac{\displaystyle\int_{D} \beta(x)\varphi_*^2}{\displaystyle\int_{D} \varphi_*^2},
\end{equation}
so that, if $\int_{D} \beta(x)\varphi_*^2>0$, we get $\sigma_1^D[-d\Delta-\beta;{\mathcal{B}}]<0$
and then, there exists positive solution $\tilde{u}_d$ for large $d$.

Moreover, multiplying (\ref{eq:classic_logistic}) by $\varphi_*^2$ and \eqref{fiestre} by $\varphi_*$, we get, rearranging terms
$$
-d\frac{\varphi_*}{\tilde{u}_d} \varphi_* \Delta \tilde{u}_d =\beta(x)\varphi_*^2-\alpha(x) \tilde{u}_d\varphi_*^2, \quad \text{and} \quad  -d\frac{\varphi_*}{\tilde{u}_d}\tilde{u}_d \Delta \varphi_* =0.
$$
Thus, adding both expression and integrating in $D$ it yields to
$$
\int_D \beta(x)\varphi_*^2-\int_D\alpha(x) \tilde{u}_d\varphi_*^2=-d\int_D \frac{\varphi_*}{\tilde{u}_d}\left(\varphi_* \Delta \tilde{u}_d - \tilde{u}_d \Delta \varphi_* \right).
$$
Applying Picone's identity, see Lemma~\ref{lemma:Picone}, with $f={\rm id}$, $u= \tilde{u}_d$ and $v=\varphi_*$ we find that
$$
-d\int_D \frac{\varphi_*}{\tilde{u}_d}\left(\varphi_*\Delta \tilde{u}_d -\tilde{u}_d\Delta  \varphi_* \right)=-d\int_D \tilde{u}_d^2 \left|\nabla\left(\frac{\varphi_*}{\tilde{u}_d}\right)\right|^2<0,
$$
and then,
\begin{equation}
\label{eq:nocero}
0<\int_D \beta(x)\varphi_*^2<\int_D\alpha(x) \tilde{u}_d\varphi_*^2.
\end{equation}
Finally, we observe that for a sufficiently large $K>0$, $K\varphi_*$ is supersolution of~(\ref{eq:classic_logistic})  if
\begin{equation}
\label{eq:k_sup}
K\alpha(x)\varphi_*\geq \beta(x).
\end{equation}
Hence, choosing the $K$ from~\eqref{eq:k_sup} it follows that
$$
\tilde{u}_d(x)\leq K\varphi_*\leq \frac{\beta_M}{\alpha_L(\varphi_*)_L}\varphi_*(x).
$$
By elliptic regularity, $\tilde{u}_d$ is bounded in $W^{2,p}(D)$ and then
$$
\tilde{u}_d\to \tilde{u}_\infty\quad\mbox{in $C^1(\overline D)$ as $d\to \infty$,}
$$
where, thanks to (\ref{eq:nocero}), $\tilde{u}_\infty>0$ and moreover it is a solution of \eqref{eq:inft_u}.
\end{proof}

\begin{remark}
\label{neu}
Note that in the particular case when $\mathcal{B}u=\partial_{{\bf n}}u$ we have that $
\tilde u_\infty=\int_D \beta(x)/\int_D \alpha(x)$, see Theorem {\rm 5.2.5} in \cite{lam-lou}.
\end{remark}

\begin{proposition}
\label{dinfilogistica.noexiste}
Problem~{\rm (\ref{eq:classic_logistic})}  does not possess positive solution for $d$ large if $\sigma_1^D[-\Delta;{\mathcal{B}}]>0$ or $\sigma_1^D[-\Delta;{\mathcal{B}}]=0$ and $\displaystyle \int_D \beta(x)\varphi_*^2<0$.
\end{proposition}

\begin{proof}
 Assume that $\sigma_1^D[-\Delta;{\mathcal{B}}]>0$. Then, by Proposition \ref{proper}
$$
\lim_{d\to\infty }\sigma_1^D[-d\Delta-\beta;{\mathcal{B}}]=+\infty,
$$
and then $\sigma_1^D[-d\Delta-\beta;{\mathcal{B}}]>0$ for a sufficiently large $d$. Hence, thanks to Proposition~\ref{prop:existence!.logistic} we get that~(\ref{eq:classic_logistic})  does not possess a positive solution for $d$ large.

On the other hand, from~\eqref{eq:noexistence.sigma0} it is straightforward to see that
$\sigma_1^D[-d\Delta-\beta;{\mathcal{B}}]>0$ if $\displaystyle\int_{D} \beta(x)\varphi_*^2<0$. Consequently, (\ref{eq:classic_logistic})  does not possess positive solution for $d$ sufficiently large.
\end{proof}

To end this section we consider the  general logistic equation with non-homogenous Robin boundary conditions
\begin{equation}
\label{spl}
\left\{\begin{array}{ll}
-d\Delta u=u(\beta(x)-\alpha(x)u) & \mbox{in $D$,}\\
\mathcal{B} u={\bf h} & \mbox{on $\p D$,}
\end{array}
\right.
\end{equation}
where  $h_i\in C(\Gamma_i)$ such that $h_i\geq 0$ on $\Gamma_i$  and $h_i\neq 0$ for some $i=1,2$.

\begin{proposition}
\label{pertur}
Problem {\rm (\ref{spl})} has a unique positive solution.
\end{proposition}

\begin{proof}
To show existence of solutions to Problem {\rm (\ref{spl})}, we use the sub-supersolution method. Take
$$
(\overline u, \underline u)=(S\psi,0)
$$
where $S>0$ is a positive constant big enough, and $\psi$ is the unique positive solution of
\begin{equation*}
\label{el}
\left\{\begin{array}{ll}
-d\D \psi+\rho \psi=1 & \mbox{in $D$,}\\
\mathcal{B} \psi=1 & \mbox{on $\partial D$.}
\end{array}
\right.
\end{equation*}
Observe that for $\rho>0$ large, we have that $\s_1^D[-d\D+\rho;\mathcal{B}]>0$. It is clear that $\underline u=0$ is subsolution (and no solution) of (\ref{spl}) and $\overline u=S\psi$ is supersolution of (\ref{spl}) provided that
$$
S\alpha(x)\psi^2\geq \psi(\beta(x)+\rho)-1 \quad x\in D,\quad S\geq \max\{\max_{\G_1}h_1,\max_{\G_2}h_2\}.
$$
This shows the existence of positive solution.

The uniqueness is straightforward because the map $u\mapsto [u(\beta(x)-\alpha(x) u)]/u$ is decreasing,
see~\cite{BrezisOswald1986}.
\end{proof}

To prove the convergence when $d\to 0$, in the spirit of Lemma 6.3 in \cite{sergiojlg} we state the following technical result.
 \begin{lemma}
\label{julimejo}
Consider $h_i\in C(\Gamma_i)$ and $j_1,j_2\in C(\overline D)$ such that $j_1(x)<j_2(x)$ for $x\in \overline D$. Then, there exists $\phi\in C^2(\overline D)$ such that
$$
j_1(x)\leq \phi(x)\leq j_2(x)\quad x\in \overline D\quad\hbox{and}\quad \mathcal{B}\phi>\mathbf{h}\quad \mbox{on $\partial D$}.
$$
\end{lemma}
 \begin{proof}
Define
$$
\delta:=\min_{x\in\overline D}(j_2(x)-j_1(x))>0.
$$
Then, there exists $\hat\phi\in C^\infty(\overline D)$ such that
$$
j_1(x)+\frac{\delta}{2}<\hat\phi(x)<j_2(x)-\frac{\delta}{2}\quad\mbox{for all $x\in \overline D$.}
$$
Also, thanks to Theorem 1.3 of \cite{sergiojlg}, there exist an open neighbourhood $\mathcal{U}\subset\mathbb{R}^N$ of $\partial D$,
a function $\psi \in C^2(\mathcal{U}; \mathbb{R})$ and a constant $\tau> 0$ such that $\psi(x) < 0$ for all $x \in  \mathcal{U}\cap D$,
$\psi(x) >0$ for all $x \in  \mathcal{U}\setminus \overline D$, $\psi(x) = 0$ for each $x \in \partial D$ and
\begin{equation}
\label{eq:volver}
\min_{\partial D}\frac{\partial \psi}{\partial {\bf n}}\geq \tau\quad\mbox{on $\partial D$.}
\end{equation}
  Take, for a constant $k\in \mathbb{R}$ to be determined later,
$$
 \hat\phi_k(x):=\hat\phi(x)-1+e^{k\psi(x)}\quad\mbox{for $x\in \mathcal{U}\cap \overline D$.}
$$
Observe that $\hat\phi_k(x)=\hat\phi(x)$ on $\partial D$, and hence, reducing $\mathcal{U}$ if necessary, we find that
$$
j_1(x)+\frac{\delta}{2}<\hat\phi_k(x)<j_2(x)-\frac{\delta}{2}\quad\mbox{for all $x\in \mathcal{U}\cap \overline D$.}
$$
On the other hand, since $\psi=0$ on $\partial D$,
$$
\begin{array}{rcl}
\mathcal{B}\hat\phi_k(x)& =&\displaystyle\mathcal{B}\hat\phi+ke^{k\psi(x)}\frac{\partial \psi}{\partial {\bf n}}+g_i (e^{k\psi(x)}-1) =\mathcal{B}\hat\phi+ke^{k\psi(x)}\frac{\partial \psi}{\partial {\bf n}} \\[10pt]
&
=&\displaystyle\mathcal{B}\hat\phi+k\frac{\partial \psi}{\partial {\bf n}}
 >h_i,\quad \hbox{on} \quad \Gamma_i,
\end{array}
$$
taking $k$ large enough by (\ref{eq:volver}). So that we finally obtain that $\mathcal{B}\hat\phi_k(x)>{\bf  h}$.
{Once we have fixed the appropriate $k$, it remains to take $\phi$ as any smooth extension of $\hat\phi_k$
from a neighborhood $\mathcal{V}$ of $\partial D$, with $\mathcal{V}\subset \mathcal{U}$, to $\overline D$ in such a way that $j_1(x)<\phi(x)<j_2(x)$ for all $x\in\overline D$.}
\end{proof}

Using the previous lemma we are able to state limiting solution when $d$ goes to zero, which turns out to be the same that the one for the homogeneous problem, see Proposition~\ref{logisnormal}.
\begin{proposition}
\label{pertur.infty}
Let  $u_d$ be the unique positive solution of {\rm (\ref{spl})}. Then
\begin{equation*}
\label{limips}
\lim_{d\to 0}u_d(x)=\frac{\beta_+(x)}{\alpha(x)}\quad\mbox{uniformly in $\overline\Omega$.}
\end{equation*}
\end{proposition}
\begin{proof}

Observe that
$$
\tilde{u}_d\leq u_d \quad\mbox{in $D$,}
$$
where $\tilde{u}_d$ is the solution of (\ref{eq:classic_logistic}). Thanks to Proposition \ref{logisnormal}, we know that
$$
\tilde{u}_d\to \frac{\beta_+(x)}{\alpha(x)}\quad \mbox{as $d\to 0$,}
$$
hence, it is enough to prove that for any $\varepsilon>0$ there exists $d_0>0$ such that for $d\in (0,d_0)$
\begin{equation}
\label{eq:aux.bound.ud}
u_d\leq  \frac{\beta_+(x)}{\alpha(x)}+\varepsilon.
\end{equation}
Now, taking
$$
j_1(x)=\frac{\beta_+(x)}{\alpha(x)}+\frac{\varepsilon}{2},\qquad j_2(x)=\frac{\beta_+(x)}{\alpha(x)}+\varepsilon
$$
and applying Lemma \ref{julimejo}, there exists $\phi\in C^2(\overline D)$ such that
$$
j_1(x)\leq \phi (x)\leq j_2(x)\quad\mbox{in $D$,}\qquad \mathcal{B}\phi\geq \mathbf{h}\quad\mbox{on $\partial D$}.
$$
Observe that
$$
\phi (x)(\beta(x)-\alpha(x)\phi (x)) =\alpha(x)\phi (x)\left(\frac{\beta(x)}{\alpha(x)}-\phi (x)\right)\leq -\frac{\varepsilon}{2}\alpha(x)\phi (x)<0.
$$
Hence,
$$
-d\Delta \phi\geq \phi (x)(\beta(x)-\alpha(x)\phi (x)) ,
$$
for $d$ sufficiently small. Consequently, since $\mathcal{B}\phi\geq {\bf h}$ on $\partial D$, we get that $\phi$ is a supersolution of~(\ref{spl}), and we deduce \eqref{eq:aux.bound.ud}.
\end{proof}

\section{Interface eigenvalue problems}
\label{Sec:Interface_eigenvalue}
In this section (and in the next one) we extend the results of the previous section to interfaces problems. In particular, given the domain configuration~\eqref{eq:domain} we first focus  on the
eigenvalue problem
\begin{equation}
\label{eq:eign_prob_membrane}
\left\{\begin{array} {l@{\quad}l}
-d\Delta \varphi_i + c_i(x)
\varphi_i=\nu \varphi_i&  \hbox{in } \Omega_i,\\
\mathcal{I}({\bf \Phi})=0 &\hbox{on } \Sigma\cup\Gamma, \end{array}\right.
\end{equation}
 where we assume that  $c_i\in C(\overline\Omega_i)$ and such that ${\bf \Phi}=(\varphi_1,\varphi_2)^T$.
Following~\cite{PAC-Bran},
to deal with  convergence of solutions, throughout this and the next sections  we consider the different functional spaces adapted to the interface problem. In particular the space
 of continuous functions with the flux condition on $\Gamma$ will be defined as
$$
\mathcal{C}_\Sigma(\Omega):=\{{\bf \Psi}\in C(\overline\Omega_1)\times C(\overline\Omega_2) : \dfrac{\p \psi_1}{\p {\bf n_1}}= \dfrac{\p \psi_2}{\p {\bf n_1}}= \gamma_i ( \psi_2- \psi_1)\  \hbox{on}\  \Sigma \quad \hbox{and}\quad  \dfrac{\partial \psi_2}{\partial {\bf n_2}}=0 \ \hbox{on}\ \Gamma\}.
$$
Analogously we define the spaces of continuously differentiable and H\"{o}lder continuous functions $\mc{C}_{\Sigma}^{1}$, $\mc{C}_{\Sigma}^{1+\eta}$,  $\dots$
for some $\eta\in (0,1]$. We also define
 $\mathcal{L}_\Sigma^p$, $\mathcal{H}_\Sigma^1$ and  $\mathcal{W}_\Sigma^{2,p}$ as  the product spaces $L^p(\Omega_1)\times L^p(\Omega_2)$,  $H^1(\Omega_1)\times H^1(\Omega_2)$ and $W^{2,p}(\Omega_1)\times W^{2,p}(\Omega_2)$, respectively, together with the boundary condition.
The norm of a function {\bf u} is defined as the sum of the norms of $u_i$
in the respective spaces; see \cite{PAC-Bran} for further details on the description of those elements.

 We state several asymptotic results when moving the diffusion parameter $d$  to zero or infinity. To this aim  we denote by
$$\Lambda_1(-d\Delta+c_1,-d\Delta+c_2)$$ the principal eigenvalue of problem (\ref{eq:eign_prob_membrane}).
With respect to the associated principal eigenfunction, that we will denote
${\bf \Phi}$ (or ${\bf \Phi}_d$) throughout the paper,  having that ${\bf \Phi}\in \mathcal{H}^1_\Sigma$.

 As in the previous section, for the sake of completeness, first we collect some properties of this principal eigenvalue.

\begin{lemma}
\label{infi}
Let us denote $\mathcal{L}^\infty=L^\infty(\Omega_1)\times L^\infty(\Omega_2)$. Then:
\begin{enumerate}
\item The map
$(d,c_1,c_2)\in (0,+\infty)\times \mathcal{L}^\infty\mapsto \Lambda_1(-d\D+c_1,-d\D+c_2)$ is continuous.
\item The map
$(c_1,c_2)\in \mathcal{L}^\infty\mapsto \Lambda_1(-d\D+c_1,-d\D+c_2)$ is increasing.
\end{enumerate}
\end{lemma}

\begin{proof}
The results are proven in~\cite{PAC-Bran, CiavolellaPaerthame, Ciavolella} and \cite{Suarezetal}.
\end{proof}

{In the spirit of Section~\ref{Sec:scalar}, we also consider the auxiliary problem
\begin{equation}
\label{eq.lineal.sin.c}
\left\{
\begin{array}{ll}
-d\Delta \varphi_{i}=\nu \varphi_i & \mbox{in $\Omega_i$,}\\
\mathcal{I}({\bf \Phi})=0 & \mbox{on $\Sigma\cup \Gamma.$}
\end{array}
\right.
\end{equation}
However, contrary to what happens in the scalar case considered in Section~\ref{Sec:scalar}, we prove now that the principal eigenvalue to Problem~\eqref{eq.lineal.sin.c} is always $0$. Hence, the dependence of the behaviour of $\Lambda_1(-d\Delta+c_1,-d\Delta+c_2)$ as $d\to\infty$ on $\Lambda_1(-\Delta, -\Delta)$ will differ from what was shown for the scalar case (see Proposition~\ref{prop.scalar.eigenvalue} and Proposition~\ref{prop.eigenvalue.d.to.infinity} below) showing the existence
of some kind of effect coming from the interface region.

\begin{lemma}
\label{lemma:membrane.laplace}
Let $\Lambda_1(-d\Delta, -d\Delta)$ and ${\bf \Phi}$ be the principal eigenvalue and eigenfunction of~\eqref{eq.lineal.sin.c} respectively. Then $\Lambda_1(-d\Delta, -d\Delta)=0$ and $\varphi_1=\varphi_2=K>0$, with $K$ constant.
\end{lemma}

\begin{proof}
 It is straightforward to see that if $\nu=0$ then ${\bf \Phi}={\bf K}=(K,K)^T$ is a solution of~\eqref{eq.lineal.sin.c}. And hence, $\nu=0$ is an eigenvalue of~\eqref{eq.lineal.sin.c} and ${\bf \Phi}$ is its corresponding eigenfunction.

 On the other hand, let ${\bf \Phi}$ be any other eigenfunction associated with the eigenvalue $0$. We divide the first equation of (\ref{eq.lineal.sin.c}) by $\varphi_{1}$ and the second one by $(\gamma_2/\gamma_1) \varphi_{2}$, integrating in $\Omega_i$, respectively, and adding both equations, we obtain
$$
\int_{\O_1}\frac{|\nabla \varphi_{1}|^2}{\varphi_{1}^2}+\frac{\g_1}{\g_2}\int_{\O_2}\frac{|\nabla \varphi_{2}|^2}{\varphi_{2}^2}+\g_1\int_\Sigma\frac{(\varphi_{2}-\varphi_{1})^2}{\varphi_{1}\varphi_{2}}
 =0,
$$
proving that ${\bf \Phi}={\bf K}$ is a solution.
\end{proof}

\begin{lemma}
  \label{lemma:bounds.Lambda1}
 Let $\Lambda_1(-d\Delta+c_1,-d\Delta+c_2)$ be the principal eigenvalue of~\eqref{eq:eign_prob_membrane}. Then:
\begin{equation}
\label{cota1}
\min\{(c_1)_L,(c_2)_L\}\leq \Lambda_1(-d\Delta+c_1,-d\Delta+c_2)\leq \max\{(c_1)_M,(c_2)_M\},
\end{equation}
and
\begin{equation}
\label{cota2}
 \Lambda_1(-d\Delta+c_1,-d\Delta+c_2)\leq \frac{\displaystyle\gamma_2\int_{\Omega_1}c_1+\gamma_1\int_{\Omega_2}c_2}{\displaystyle\g_2|\Omega_1|+\gamma_1|\Omega_2|}.
\end{equation}
\end{lemma}

\begin{proof}
  Observe that,  since $\Lambda_1$ is increasing in $c_1$ and $c_2$ we have,
$$
\begin{aligned}
   \Lambda_1(-d\Delta+c_1,-d\Delta+c_2)&\leq  \Lambda_1(-d\Delta+(c_1)_M,-d\Delta+(c_2)_M)
   \\ &\leq \Lambda_1(-d\Delta,-d\Delta)+\max\{(c_1)_M,(c_2)_M\}=\max\{(c_1)_M,(c_2)_M\},
\end{aligned}$$
where we have used that  $ \Lambda_1(-d\Delta,-d\Delta)=0$ .
Analogously, we find
$$
 \Lambda_1(-d\Delta+c_1,-d\Delta+c_2)\geq \min\{(c_1)_L,(c_2)_L\},
$$
and one can deduce (\ref{cota1}).

On the other hand from~\eqref{eq:eign_prob_membrane} we have that
$$
-d\frac{\Delta\varphi_i}{\varphi_i}+c_i=\Lambda_1(-d\Delta+c_1,-d\Delta+c_2)\quad\mbox{in $\Omega_i$.}
$$
Integration by parts in $\Omega_i$ yields
$$
 d\left(-\int_{\Omega_i}\frac{|\nabla \varphi_i|^2}{\varphi_i^2}-\int_\Sigma\frac{\partial_{{\bf n}_i}\varphi_i}{\varphi_i}\right)+\int_{\Omega_i}c_i=\Lambda_1(-d\Delta+c_1,-d\Delta+c_2)|\Omega_i|.
$$
Now, we add the equation for $i=1$ and the equation fo $i=2$ multiplied by $\gamma_1/\gamma_2$ and we obtain
\begin{align*}
& -d\left(\int_{\Omega_1}\frac{|\nabla \varphi_1|^2}{\varphi_1^2}+\frac{\gamma_1}{\gamma_2}\int_{\Omega_2}\frac{|\nabla \varphi_2|^2}{\varphi_2^2}+\gamma_1\int_\Sigma\frac{(\varphi_2-\varphi_1)^2}{\varphi_1\varphi_2}\right)
+\int_{\Omega_1}c_1+\frac{\gamma_1}{\gamma_2}\int_{\Omega_2}c_2 \\ & = \Lambda_1(-d\Delta+c_1,-d\Delta+c_2)\left(|\Omega_1|+\frac{\gamma_1}{\gamma_2}|\Omega_2|\right),
\end{align*}
whence we deduce (\ref{cota2}).
\end{proof}

To study the behaviour of the first eigenvalue when $d\to 0$, we introduce
 the  homogeneous problem associated to~\eqref{eq:eign_prob_membrane},
\begin{equation}
\label{eq:homogeneous_prob_membrane}
\left\{\begin{array} {l@{\quad}l}
-d\Delta u_i + c_i(x)
u_i=0,&  \hbox{in } \Omega_i,\\
\mathcal{I}({\bf u})=0, &\hbox{on } \Sigma\cup\Gamma, \end{array}\right.
\end{equation} and state the existence of solutions,  the analogous of Proposition~\ref{prop:existence!.logistic} for this interface problem.
To do so, we define:

\begin{definition}
\label{super}
Given $\overline{\bf u}\geq 0$ in $\Omega$, {$\overline{\bf u}\in \mathcal{W}^{2,p}$}, $p>N$, we say that $\overline{\bf u}$ is a strict supersolution of~\eqref{eq:homogeneous_prob_membrane} if
$$
-d\Delta\overline {\bf u} +{\bf c}(x)\overline {\bf u} \geq 0\quad \mbox{in $\Omega$,}\qquad \left\{\begin{array}{ll}
\partial_{\bf n_1}u_1\geq\gamma_1(u_2-u_1)&  \hbox{on } \Sigma,\\
\partial_{\bf n_1}u_2\leq\gamma_2(u_2-u_1) & \hbox{on } \Sigma,\\
\partial_{\bf n_2}u_2\geq0&  \hbox{on } \Gamma,
\end{array}
\right.
$$
and some of these inequalities are strict.
\end{definition}

\begin{lemma}
\label{lem:existence.interface.eingenvale}
There exists a positive  strict supersolution $\overline{\bf u}$ of~\eqref{eq:homogeneous_prob_membrane} if and only if $$ \Lambda_1(-d\D+c_1,-d\D+c_2)>0.$$
\end{lemma}

\begin{proof}
The proof can be read in~\cite{Suarezetal}.
\end{proof}

 Once we have stated the main properties of the principal eigenvalue, we deal now with its behaviour for small and large diffusion.
\begin{proposition}
\label{dauto}
Let $\Lambda_1(-d\Delta+c_1,-d\Delta+c_2)$ be the principal eigenvalue of~\eqref{eq:eign_prob_membrane}. Then:
\begin{equation}
\label{dacero}
\lim_{d\to 0} \Lambda_1(-d\Delta+c_1,-d\Delta+c_2)=\min\{(c_1)_L,(c_2)_L\}.
\end{equation}
\end{proposition}

\begin{proof}
Without loss of generality, we might assume that
$$(c_1)_L<(c_2)_L.
$$
By (\ref{cota1}) we get that
$$
\Lambda_1(-d\Delta+c_1,-d\Delta+c_2)\geq (c_1)_L.
$$
Reasoning by contradiction we assume that there exists $\varepsilon>0$ such that
$$
{\lim_{d\to 0}}\Lambda_1(-d\Delta+c_1,-d\Delta+c_2)>(c_1)_L+\varepsilon,
$$
and then, for a sequence $\{d_n\}$ converging to zero, we find that
$$
\Lambda_1(-d_n\Delta+c_1-(c_1)_L-\varepsilon,-d_n\Delta+c_2-(c_1)_L-\varepsilon)>0.
$$
This implies, due to Lemma~\ref{lem:existence.interface.eingenvale}, that there exists a {strict positive} function ${\bf \Phi}_n$, such that
$$
\left\{
\begin{array}{ll}
(-d_n\Delta+c_1-(c_1)_L-\varepsilon)\varphi_{1,n}\geq0&\quad\mbox{in $\Omega_1$,}\\
(\partial_{{\bf n}_1}+\gamma_1)\varphi_{1,n}\geq \gamma_1\varphi_{2,n}>0&\quad\mbox{on $\Sigma$.}
\end{array}
\right.
$$
On the other hand, for $\varepsilon >0$,  there is a ball $B$ such that
$$
c_1(x)<(c_1)_L+\frac{\varepsilon}{2},\qquad x\in B.
$$
Then,
$$
\left(-d_n\Delta-\frac{\varepsilon}{2}\right)\varphi_{1,n}>(-d_n\Delta+c_1-(c_1)_L-\varepsilon)\varphi_{1,n}\geq 0\quad\mbox{in $B$,} \qquad \mbox{and}\qquad
\varphi_{1,n}>0\quad\mbox{on $\p B$.}
$$
Consequently, denoting $\sigma_1^B[-d_n\D-\frac{\varepsilon}{2};\mathcal{D}]$ by the principal eigenvalue under Dirichlet boundary conditions, it follows that
$$
0<\sigma_1^B[-d_n\Delta-\frac{\varepsilon}{2};\mathcal{D}]=\sigma_1^B[-d_n\Delta;\mathcal{D}]-\frac{\varepsilon}{2},
$$
which is a contradiction, since $\displaystyle\lim_{d_n\to 0}\sigma_1^B[-d_n\Delta;\mathcal{D}]=0.$
\end{proof}

We conclude this section by stating the behaviour of the principal eigenvalue when the diffusion coefficient tends to infinity. As mentioned above, the result is independent of the sign of $\Lambda_1(-d\Delta,-d\Delta)$  and of the associated eigenfunctions.

\begin{proposition}
\label{prop.eigenvalue.d.to.infinity} Let $\Lambda_1(-d\Delta+c_1,-d\Delta+c_2)$ be the principal eigenvalue of~\eqref{eq:eign_prob_membrane}. Then
  \begin{equation}
\label{dainfinito}
\lim_{d\to \infty} \Lambda_1(-d\Delta+c_1,-d\Delta+c_2)= \frac{\displaystyle\gamma_2\int_{\Omega_1}c_1+\gamma_1\int_{\Omega_2}c_2}{\gamma_2|\Omega_1|+\gamma_1|\Omega_2|}.
\end{equation}
\end{proposition}

\begin{proof}
Take a sequence of positive eigenfunctions ${\bf \Phi}_d$,
associated with the principal eigenvalue $\Lambda_1(-d\Delta+c_1,-d\Delta+c_2)$ such that
$$
 \|\varphi_{i,d}\|_\infty=1.
$$

Since $\varphi_{i,d}$ is bounded in $L^\infty(\O_i)$, then it is bounded in $L^p(\O_i)$ for all $p\geq 1$. Hence, by elliptic regularity,
$\varphi_{i,d}$ is bounded in $W^{1,p}(\O_i)$ for $1<p<\infty$, so that $\varphi_{i,d}$ is bounded in $W^{1-1/p,p}(\partial\Omega_i)$. Then, we conclude that $\varphi_{i,d}$ is bounded in $W^{2,p}(\O_i)$. Also,
we might conclude that $\varphi_{i,d}$ is bounded in $C^{1,\alpha}(\overline\Omega_i)$; see \cite{GT} for further details. Going back to the equations~\eqref{eq:eign_prob_membrane}, $\varphi_{i,d}$ is bounded in $C^{2,\alpha}(\overline\O_i)$. Then,  we can pass to the limit as $d\to\infty$ and get
\begin{equation*}
\label{eq:lim_infinito}
{\bf \Phi}_d\to {\bf \Phi}_\infty\quad\mbox{in $\mathcal{W}^{2,p}_\Sigma$ and in $\mathcal{C}^1_\Sigma$}
\end{equation*}
and dividing each equation by $d$ we have that  ${\bf \Phi}_\infty$ is the principle eigenfunction of~\eqref{eq.lineal.sin.c}, with $d=1$. Hence, from Lemma~\ref{lemma:membrane.laplace}, we have  $\varphi_{i,\infty}=1$.

On the other hand, multiplying  the second equation of \eqref{eq:eign_prob_membrane} by $\gamma_1/\gamma_2$, integrating in $\Omega_i$, respectively, and adding both equalities we arrive at
\begin{equation*}
\label{paso}
\begin{aligned}
-d\int_{\Omega_1}\Delta \varphi_{1,d} -d\frac{\gamma_1}{\gamma_2}\int_{\Omega_2}\Delta \varphi_{2,d} &+
\int_{\Omega_1}\varphi_{1,d} c_1+\frac{\gamma_1}{\gamma_2}\int_{\Omega_2}\varphi_{2,d} c_2\\&= \Lambda_1(-d\Delta+c_1,-d\Delta+c_2)\left(\int_{\Omega_1}\varphi_{1,d}+\frac{\gamma_1}{\gamma_2}\int_{\Omega_2}\varphi_{2,d} \right).
\end{aligned}
\end{equation*}
It is easy to see, using Gauss Divergence Theorem, that
\begin{equation}
  \label{eq:divergence}
\int_{\Omega_1}\Delta \varphi_{1,d} +\frac{\gamma_1}{\gamma_2}\int_{\Omega_2}\Delta \varphi_{2,d}=0.
\end{equation}
Therefore,
$$\displaystyle
\Lambda_1(-d\Delta+c_1,-d\Delta+c_2)=\dfrac{\displaystyle \int_{\Omega_1}\varphi_{1,d} c_1+\dfrac{\gamma_1}{\gamma_2}\int_{\Omega_2}\varphi_{2,d} c_2}{\displaystyle\left(\int_{\Omega_1}\varphi_{1,d}+\frac{\gamma_1}{\gamma_2}\int_{\Omega_2}\varphi_{2,d} \right)}.
$$
Finally, passing to the limit as $d$ goes to infinity we prove~(\ref{dainfinito}).
\end{proof}

\section{Interface logistic equation}
\label{Sec:interface_logistic}

Consider  now the general interface logistic problem~\eqref{eq:general-system}--\eqref{eq:neumann-condition}. That is
  \begin{equation}
\label{eqlogisdinfi} \left\{\begin{array} {l@{\quad}l} -d\Delta
u_{i} =u_{i}(\beta_i(x)-\alpha_i (x)u_{i})& \mbox{in $\Omega_i$},\\
\mathcal{I}({\bf u})=0 & \mbox{on $\Sigma\cup \Gamma$},
\end{array}\right.
\end{equation}
 where $\alpha_i,\beta_i\in C(\overline\Omega_i) $, with $\alpha_i(x)>0$ for all $x\in\overline\Omega_i$ and $(\beta_i)_M>0$.

Thus, we show the existence and uniqueness of positive solutions for the interface nonlinear problem (\ref{eqlogisdinfi}), which can be seen after applying the appropriate existence results shown in \cite{PAC-Bran} and \cite{Suarezetal}.
\begin{proposition}
\label{exisuni}
There exists a positive solution of {\rm (\ref{eqlogisdinfi})} if and only if
$$
\Lambda_1(-d\Delta-\beta_1,-d\Delta-\beta_2)<0.
$$
When the positive solution exists, it is  unique.
\end{proposition}

\begin{proof}
Let $\Lambda_1(-d\Delta-\beta_1,-d\Delta-\beta_2)<0$.  To prove this existence result we will apply the method of sub and supersolutions.

First we look for a subsolution for~{\rm(\ref{eqlogisdinfi})}. To do so, consider $\epsilon {\bf \Phi}$, for $0<\epsilon\ll 1$, and where ${\bf \Phi}$
is the principal (normalized) eigenfunction associated with the principal eigenvalue {$\Lambda_1(-d\Delta-\beta_1,-d\Delta-\beta_2)$} for problem~\eqref{eq:eign_prob_membrane}. Indeed,
$$
-d\Delta(\epsilon\varphi_i) =\epsilon\varphi_i\left(\Lambda_1(-d\Delta-\beta_1,-d\Delta-\beta_2) +\beta_i\right)< \epsilon\varphi_{i}(\beta_i-\alpha_i \epsilon\varphi_{i})
\quad \hbox{in}\quad \Omega_i,
$$
and hence,
$$\Lambda_1(-d\Delta-\beta_1,-d\Delta-\beta_2) < -\alpha_i \epsilon \varphi_{i},$$
which is satisfied for a sufficiently small  $\epsilon>0$.
On the boundary we have the equality.

On the other hand, for a sufficiently big constant $S$ we claim that $(S,S)^T$ is a supersolution to~\eqref{eqlogisdinfi}. Indeed,  substituting $S$ into the equations of \eqref{eqlogisdinfi} we obtain that if
\begin{equation}
\label{reina}
S\geq \max_{i=1,2}\{(\beta_i)_M/(\alpha_i)_L\},
\end{equation}
then
$$
0> \beta_i S- \alpha_i S^2\quad \hbox{in}\quad \Omega_i.
$$

Now, assume that there exists a positive solution $(u_1,u_2)$ of (\ref{eqlogisdinfi}), then
$$
0=\Lambda_1(-d\Delta-\beta_1+\alpha_1 u_1,-d\Delta-\beta_2+\alpha_2 u_2)>\Lambda_1(-d\Delta-\beta_1,-d\Delta-\beta_2).
$$
The uniqueness follows by \cite{Suarezetal}.
\end{proof}

\begin{proposition}
\label{prop.cota.uid}
Let  ${\bf u}_d$ be the unique positive solution of {\rm (\ref{eqlogisdinfi})}.
Then
\begin{equation}
\label{cotau}
u_{i,d}\leq\max_{i=1,2}\left\{\frac{(\beta_i)_M}{(\alpha_i)_L}\right\}.
\end{equation}
\end{proposition}
\begin{proof}
The result is straightforward, just noting that $S=\max_{i=1,2}\{(\beta_i)_M/(\alpha_i)_L$ is a supersolution of  (\ref{eqlogisdinfi}), see (\ref{reina}), and the uniqueness of positive solution of (\ref{eqlogisdinfi}).
\end{proof}

 Now we study the singular case $d\to 0$. 
\begin{theorem}
\label{dcero}
Let ${\bf u}_d$ be the unique positive solution of {\rm (\ref{eqlogisdinfi})}. Then
\begin{equation*}
\label{dtendiendocero}
\lim_{d\to 0}u_{i,d}=\frac{(\beta_i(x))_+}{\alpha_i(x)}\quad\mbox{uniformly in $\overline\Omega_i$}.
\end{equation*}
\end{theorem}

\begin{proof}
First, observe that from (\ref{dacero}) in Proposition~\ref{dauto} we have
$$
\Lambda_1(-d\Delta-\beta_1,-d\Delta -\beta_2)\to \min\{-(\beta_1)_M,-(\beta_2)_M\}<0, \qquad \hbox{as}\quad d\to 0.
$$
Hence, by Proposition~\ref{exisuni}, there exists a positive solution of~\eqref{eqlogisdinfi} for $d$ small.

Now, consider  the following uncoupled system
\begin{equation}
\label{w1-w2} \left\{\begin{array} {ll}
-d\Delta w_{i} =w_{i}(\beta_i(x)-\alpha_i (x)w_{i})& \hbox{in $\Omega_i$},\\
\partial_{{\bf n}_i} w_i+\gamma_i w_i=0 & \hbox{on $\Sigma$,}\\
\partial_{{\bf n}_2} w_2=0 & \hbox{on $\Gamma$.}
\end{array}
\right.
\end{equation}
Note that, thanks to {Proposition~\ref{prop.scalar.limit.0}}, we find that
$$
\lim_{d\to 0} \sigma_1^{\Omega_1}[-d\Delta-\beta_1;\mathcal{N}+\gamma_1]=(-\beta_1)_L=-(\beta_1)_M<0,
$$
and
$$
\lim_{d\to 0} \sigma_1^{\Omega_2}[-d\Delta-\beta_2;\mathcal{N}+\gamma_2;\mathcal{N}]=(-\beta_2)_L=-(\beta_2)_M<0.
$$
Hence, due to {Proposition~\ref{prop:existence!.logistic}}, $w_{1,d}$ and $w_{2,d}$ are the unique solutions of (\ref{w1-w2}).

On the other hand, let us consider as well the following non-homogeneous version of~\eqref{w1-w2} in which the boundary data on $\Sigma$ are replaced by
$$
\partial_{{\bf n}_1} W_1+\gamma_1 W_1=\g_2 K,\quad \mbox{and} \quad
\partial_{{\bf n}_2} W_2+\gamma_2 W_2=\gamma_1 K,
$$
where
$K=\max\{ \frac{(\beta_1)_M}{(\alpha_1)_L}, \frac{(\beta_2)_M}{(\alpha_2)_L}\} $.
Thanks to {Proposition \ref{pertur}}, these non-homogeneous equations possess a unique positive solution, $W_{1,d}$ and $W_{2,d}$, respectively. Moreover,
due to the bound (\ref{cotau}) we have that
$$
w_{i,d}\leq u_{i,d}\leq W_{i,d}\quad\mbox{in $\O_i$.}
$$
Consequently, using {Proposition \ref{logisnormal} and Proposition \ref{pertur.infty}} we complete the proof of the uniform convergence of $u_{i,d}$ as the diffusion coefficient $d$ goes to zero.
\end{proof}

 In the next result, we study the behaviour of the solution of (\ref{eqlogisdinfi}) as $d\to\infty$. As we have mentioned in Section~\ref{Sec:Interface_eigenvalue}, the fact that the principal eigenvalue of~\eqref{eq.lineal.sin.c} is $\Lambda_1(-d\Delta,-d\Delta)=0$, produces a different behaviour between the membrane solution and the scalar one; see Section~\ref{Sec:scalar}. In fact, even when
 the next result is similar to Propositions~\ref{dinfilogistica} and~\ref{dinfilogistica.noexiste}, here the condition does not depend on the eigenvalue as the previous dependence on
 $\sigma^D_1[-\Delta;\mathcal{B}]$. Moreover, the limit behaviour is a number, whereas for the scalar problem it is a function depending on the principal eigenfunctions.

\begin{theorem}
  \label{otrocerogene}
Let  ${\bf u}_d$ be the unique positive solution of {\rm (\ref{eqlogisdinfi})}. If
\begin{equation}
\label{hipo}
\gamma_2\int_{\Omega_1}\beta_1+\gamma_1\int_{\Omega_2}\beta_2>0.
\end{equation}
Then,
\begin{equation}
\label{comportamientodinfty}
\lim_{d\to \infty}u_{i,d}=\frac{\displaystyle \gamma_2 \int_{\Omega_1}\beta_1 +\gamma_1\int_{\Omega_2}\beta_2}{\displaystyle \gamma_2 \int_{\Omega_1}\alpha_1+\gamma_1\int_{\Omega_2}\alpha_2}\quad\mbox{in $C(\overline\Omega_i)$.}
\end{equation}
\end{theorem}

\begin{proof}
Observe that, thanks to Proposition~\ref{prop.eigenvalue.d.to.infinity},  it follows that
$$
\lim_{d\to\infty}\Lambda_1(-d\Delta-\beta_1,-d\Delta-\beta_2)= -\frac{\displaystyle\gamma_2\int_{\Omega_1}\beta_1+\gamma_1\int_{\Omega_2}\beta_2}{\gamma_2|\Omega_1|+\gamma_1|\Omega_2|}.
$$
Hence,  assuming that the inequality (\ref{hipo}) is satisfied, one has,  that $\Lambda_1(-d\Delta-\beta_1,-d\Delta-\beta_2)<0$ for $d$ sufficiently  and from Proposition~\ref{exisuni}, there exists a positive solution ${\bf u}_i$ of {\rm (\ref{eqlogisdinfi})}.

Now, we present two different proofs to characterize that limit solution \eqref{comportamientodinfty}:

\noindent{\it Proof 1.} Due to estimation (\ref{cotau}) it follows that $u_{i,d}$ is bounded in $L^\infty(\Omega_i)$. Thus, the right hand sides of the equations
\begin{equation*}
\label{eq:funtion_i}
u_{i,d}(\beta_i(x)-\alpha_i(x)u_{i,d})
\end{equation*}
are bounded in $L^p(\Omega_i)$ for all $p\geq 1$. Consequently, using the same elliptic regularity argument  as performed in the proof of Proposition~\ref{prop.eigenvalue.d.to.infinity} (see \cite{GT} for further details) we conclude that
$u_{i,d}$ are bounded  in $C^{2,\alpha}(\overline\Omega_i)$. Therefore
$$
u_{i,d}\to u_{i,\infty}\quad\mbox{in $C^2(\overline\Omega_i)$}.
$$
Moreover, since $d\to \infty$, we conclude that $u_{i,\infty}$  is a solution of~\eqref{eq.lineal.sin.c} with $\nu=0$ and hence (see Lemma~\ref{lemma:membrane.laplace}),
$$
u_{1,\infty}=u_{2,\infty}=L.
$$
Arguing again as in the proof of Proposition~\ref{prop.eigenvalue.d.to.infinity} we take the first equation of problem (\ref{eqlogisdinfi}), integrate in $\Omega_1$, and also multiply the second equation by $\gamma_1/\gamma_2$, integrate in $\Omega_2$, and add both equalities to get
$$
L=\displaystyle\frac{\displaystyle\gamma_2\int_{\Omega_1}\beta_1+\gamma_1\int_{\Omega_1}\beta_2}
{\displaystyle\gamma_2\int_{\Omega_1}\alpha_1+\gamma_1\int_{\Omega_1}\alpha_2}.
$$

\noindent{\it Proof 2.} The second proof we present here has the particular characteristic of not needing regularity for the functions $\alpha_i$ and $\beta_i$ and it in the spirit of Theorem 5.2.5 in \cite{lam-lou}. Thus, we multiply the first equation of problem (\ref{eqlogisdinfi}) by $u_{1,d}$ and the second one by $(\gamma_1/\gamma_2)u_{2,d}$, integrate in $\Omega_i$ respectively. Adding the two expressions yields
\begin{align*}
 d&\left( \int_{\Omega_1}|\nabla u_{1,d} |^2   +(\gamma_1/\gamma_2)\int_{\Omega_2}|\nabla u_{2,d} |^2 +(\gamma_1/\gamma_2)\int_{\Sigma}(u_{2,d}-u_{1,d})^2\right)
\\   &= \int_{\Omega_1}f_1(x,u_{1,d})u_{1,d}+(\gamma_1/\gamma_2)\int_{\Omega_2}f_2(x,u_{2,d})u_{2,d},
\end{align*}
where $f_i(x,u_{i,d})=u_{i,d}(\beta_i(x)-\alpha_i(x)u_{i,d})$.
 Thus, thanks to the {bound (\ref{cotau})} we have that
$$
\int_{\O_i}|\nabla u_{i,d} |^2\leq \frac{C}{d}\qquad\hbox{and}\qquad  \int_{\Sigma}(u_{2,d}-u_{1,d})^2\leq \frac{C}{d}.
$$
Also, due to Poincar\'e's inequality we find that
$$
\int_{\O_i}\left|u_{i,d}-\frac{1}{|\O_i|}\int_{\O_i}u_{i,d}\right|^2\leq \int_{\O_i}|\nabla u_{i,d} |^2\leq \frac{C}{d}.
$$
Therefore,
denoting
$$
M_i:= \frac{1}{|\Omega_i|}\int_{\Omega_i}u_{i,d},
$$
it follows that
$$
u_{i,d}-M_i\to 0\qquad\mbox{in $L^2(\Omega_i)$, as $d\to \infty$.}
$$
As a consequence, we arrive at
\begin{equation}
\label{eq:H1_limit}
u_{i,d}-M_i\to 0\qquad\mbox{in $H^1(\Omega_i)$, as $d\to \infty$,}
\end{equation}
and thanks to the trace inequality we have as well that
$$
u_{i,d}-M_i\to 0\qquad\mbox{in $L^2(\Sigma)$, as $d\to \infty$.}
$$
We can conclude that
\begin{equation}
\label{esen}
M_1-M_2\to 0,\quad \text{and} \quad u_i-M_i\to 0 \qquad \mbox{as $d\to \infty$.}
\end{equation}
Furthermore, recall that, as in~\eqref{eq:divergence},
$$
0=\int_{\Omega_1}f_1(x,u_{1,d}) +\frac{\gamma_1}{\gamma_2}\int_{\Omega_2}f_2(x,u_{2,d}).
$$
Hence, we have
$$\begin{aligned}
-\int_{\Omega_1}\left(u_{1,d}-M_1\right)(\beta_1-\alpha_1 u_{1,d})&-\frac{\gamma_1}{\gamma_2}\int_{\Omega_2}\left(u_{2,d}-M_2\right)(\beta_2-\alpha_2 u_{2,d})\\
&=M_1\left[\int_{\Omega_1}\beta_1-\int_{\Omega_1}\alpha_1 u_{1,d}\right]+\frac{\gamma_1}{\gamma_2}M_2\left[\int_{\Omega_2}\beta_2-\int_{\Omega_2}\alpha_2 u_{2,d}\right].
\end{aligned}
$$
Then, thanks to the limit \eqref{eq:H1_limit} it follows that
$$
\lim_{d\to \infty}\left\{ M_1\left[\int_{\Omega_1}\beta_1-\int_{\Omega_1}\alpha_1 u_{1,d}\right]+\frac{\gamma_1}{\gamma_2}
M_2\left[\int_{\Omega_2}\beta_2-\int_{\Omega_2}\alpha_2 u_{2,d}\right]\right\}=0.
$$
Now, using (\ref{esen}) we get that
\begin{equation}
\label{eq:proof2.aux}
\lim_{d\to \infty}\left(\gamma_2\int_{\Omega_1}\alpha_1 u_{1,d}+\gamma_1 \int_{\Omega_2}\alpha_2u_{2,d}\right)=\gamma_2\int_{\Omega_1}\beta_1+\gamma_1\int_{\Omega_2}\beta_2.
\end{equation}
On the other hand,
$$
\begin{aligned} \g_2\int_{\Omega_1}\alpha_1&\left(u_{1,d}-M_1\right)+\gamma_1\int_{\Omega_2}\alpha_2\left(u_{2,d}-M_2\right)=
\\
&\gamma_2\int_{\Omega_1}\alpha_1 u_{1,d}+\gamma_1 \int_{\Omega_2}\alpha_2u_{2,d}-\gamma_2\int_{\Omega_1}\alpha_1\left[ M_1-M_2\right]-M_2\left[\gamma_2\int_{\Omega_1}\alpha_1+\gamma_1\int_{\Omega_2}\alpha_2\right].
\end{aligned}
$$
Passing to the limit and using~\eqref{eq:proof2.aux}, we get
$$
M_2:=\frac{1}{|\Omega_2|}\int_{\Omega_2}u_{2,d}\to \frac{\displaystyle\gamma_2\int_{\Omega_1}\beta_1+\gamma_1\int_{\Omega_2}\beta_2}{\displaystyle\gamma_2\int_{\Omega_1}\alpha_1+\gamma_1\int_{\Omega_2}\alpha_2},
$$
from where we conclude~\eqref{comportamientodinfty} using again~\eqref{esen}.
\end{proof}

\begin{lemma}
  Problem {\rm(\ref{eqlogisdinfi})} does not possess positive solution  when $d$ tends to infinity if
$$
\gamma_2\int_{\Omega_1}\beta_1+\gamma_1\int_{\Omega_2}\beta_2<0.
$$
\end{lemma}
\begin{proof}
  It is a direct consequence of Proposition~{\rm\ref{prop.eigenvalue.d.to.infinity}} and Proposition{\rm~\ref{exisuni}}.
\end{proof}

\section{Interface logistic equation with constant growth rates}
\label{Sec:interface_growth_rates}

In this section we consider again a logistic interface problem~\eqref{eq:general-system}--\eqref{eq:neumann-condition}, however assuming constant growth rates $\beta_i(x)\equiv \lambda_i\in \mathbb{R}$:
\begin{equation}
\label{eq:membrane.logis}
\left\{\begin{array} {l@{\quad}l}
-\Delta u_{i} =u_{i}(\lambda_i-\alpha_i (x)u_{i})& \hbox{in $\Omega_i$},\\
\mathcal{I}({\bf u})=0 & \hbox{on $\Sigma\cup  \Gamma$.}
\end{array}
\right.
\end{equation}
Since $d=1$, we simplify the notation introduced in the previous sections and denote the principal eigenvalue of~\eqref{eq:eign_prob_membrane} as
$$
\Lambda_1(c_1,c_2):=\Lambda_1(-\Delta+c_1,-\Delta+c_2).
$$

Next proposition is a special case of  Proposition~\ref{exisuni} with $d=1$ and $\beta_i=\lambda_i$.
\begin{proposition}
\label{exislogismem}
There exists a positive solution of~{\rm(\ref{eq:membrane.logis})} if and only if
\begin{equation}
\label{condinecesufi}
\Lambda_1(-\lambda_1,-\lambda_2)<0.
\end{equation}
When the solution exists, it is  unique.
\end{proposition}

In what follows, let us consider equations {\rm(\ref{w1-w2})} with $d=1$ and $\beta_i(x)=\lambda_i$. that is
\begin{equation}
\label{w1-w2.sec5} \left\{\begin{array} {ll}
-\Delta u_{i} =u_{i}(\lambda_i-\alpha_i (x)u_{i})& \hbox{in $\Omega_i$},\\
\partial_{{\bf n}_i} u_i+\gamma_i u_i=0 & \hbox{on $\Sigma$,}\\
\partial_{{\bf n}_2} u_2=0 & \hbox{on $\Gamma$,}
\end{array}
\right.
\end{equation}
and denote $u_{\lambda_i}$ its  unique positive solution.

\begin{corollary}
\label{prop.cota.lambda}
Let ${\bf \Theta}_\lambda$ be the unique solution of~{\rm\eqref{eq:membrane.logis}}. Then
\begin{equation}
\label{cotaimpor}
u_{\lambda_i}\leq \theta_{\lambda_i}\leq \max\left\{\frac{\lambda_1}{(\alpha_1)_L},\frac{\lambda_2}{(\alpha_2)_L}\right\}\quad\mbox{in $\Omega_i$.}
\end{equation}
\end{corollary}
\begin{proof}
The result is straightforward, just by observing that $\theta_{\l_i}$ is a supersolution of equations~(\ref{w1-w2.sec5}), for $i=1,2$, respectively.
\end{proof}
In the following result we characterize the set of values $(\lambda_1,\lambda_2)\in \mathbb{R}^2$ verifying condition (\ref{condinecesufi}), i.e. values of $(\lambda_1,\lambda_2)$
providing us with the existence of positive solutions.
To this end, we first study the equation
 \begin{equation}
\label{ecuimpor}
\Lambda_1(-\lambda_1,-\lambda_2)=0.
\end{equation}
Equation (\ref{ecuimpor}) defines a curve in $\mathbb{R}^2$ that we now describe, denoting
$$
\sigma_1:= \sigma_1^{\Omega_1}[-\Delta;\mathcal{N}+\gamma_1],\qquad \sigma_2:= \sigma_1^{\Omega_2}[-\Delta;\mathcal{N}+\gamma_2;\mathcal{N}],
$$
where $\sigma_1$ and $\sigma_2$ stand for principal eigenvalues of two uncoupled problems as described in~\eqref{eq:princ_eign}.
Observe that $\sigma_1,\sigma_2>0$, (see Figure \ref{fig1}  and \cite{Suarezetal2}).
\begin{figure}[ht!]
\includegraphics[scale=0.4]{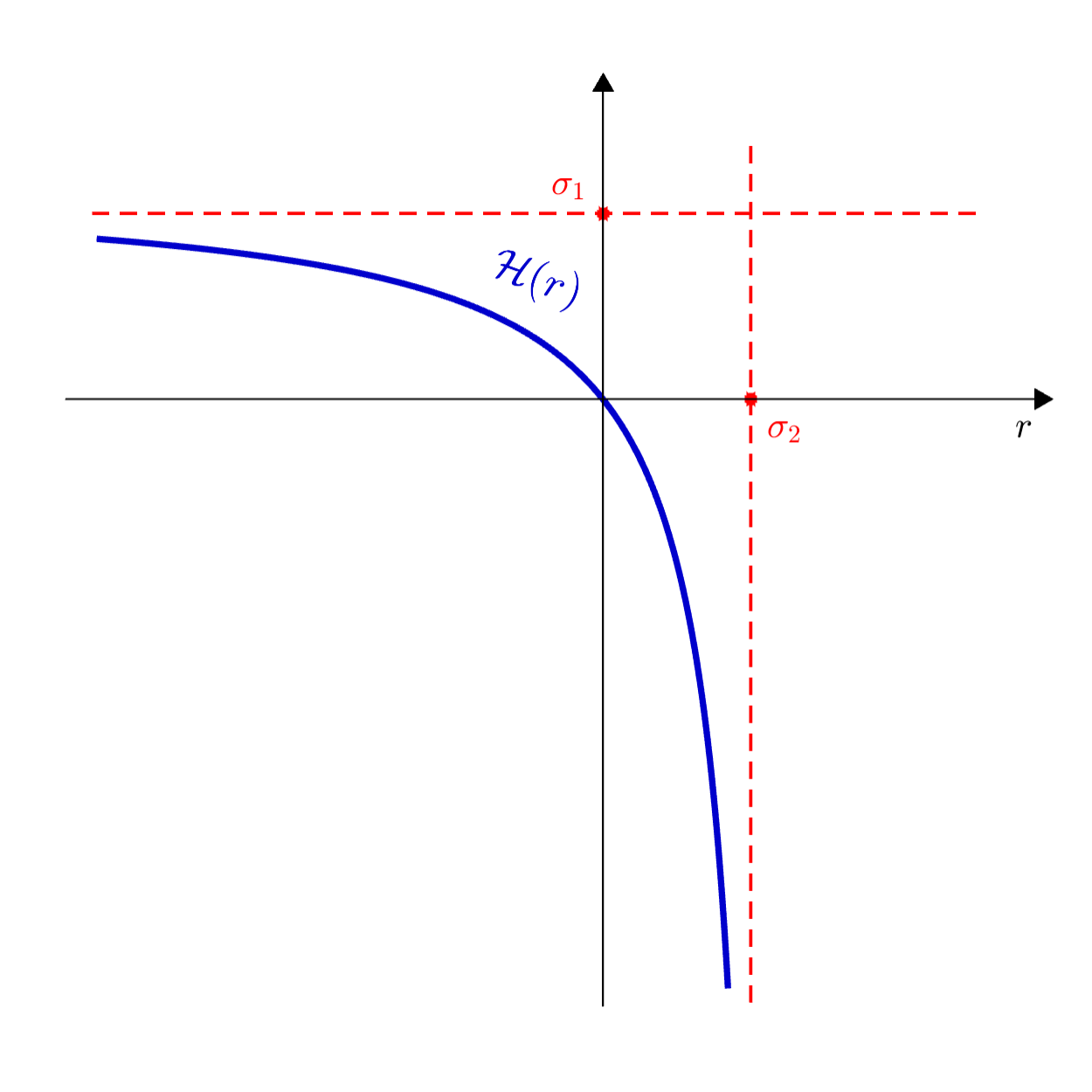}
\caption{Curve determining the eigenvalues of (\ref{eq:membrane.logis}).}
\label{fig1}
\end{figure}

\begin{lemma}
\label{autovalor}
There exists a continuous and decreasing map $
\mathcal{H}:(-\infty,\sigma_2)\mapsto \mathbb{R},
$ with
$$
\lim_{r\to -\infty}\mathcal{H}(r)=\sigma_1,\qquad \lim_{r\to \sigma_2}\mathcal{H}(r)=-\infty,
$$
such that:
\begin{enumerate}
\item $\Lambda_1(-\nu_1,-\nu_2)=0$ if and only if $\nu_1=\mathcal{H}(\nu_2)$.
\item If $\nu_2\geq \sigma_2$, then $\Lambda_1(-\nu_1,-\nu_2)<0$ for all $\nu_1\in \mathbb{R}$.
\item If $\nu_2< \sigma_2$ then
$
\Lambda_1(-\nu_1,-\nu_2)<0 \ \mbox{(resp. $>0$)\quad  if\quad } \nu_1> \mathcal{H}(\nu_2)\  \mbox{(resp. $\nu_1< \mathcal{H}(\nu_2)$.)}
$
\end{enumerate}
\end{lemma}

\begin{proof} The proof of this result can be found in \cite{Suarezetal2}.
\end{proof}

\begin{remark}
\label{dis}
Observe that using Lemma~{\rm\ref{autovalor}}  we can conclude existence and uniqueness of solutions to~\eqref{eq:membrane.logis}, see Proposition~{\rm\ref{exislogismem}} as follows:
\begin{enumerate}
\item If $\l_2\geq \sigma_2$, {\rm(\ref{condinecesufi})} holds for all $\l_1\in\mathbb{R}$.
\item If $\l_2< \sigma_2$, {\rm(\ref{condinecesufi})} holds $\l_1>\mathcal{H}(\l_2)$.
\end{enumerate}
\end{remark}

In the rest of this section we are going to analyse the behaviour of the solutions depending on the growth rates, assuming either that the growth rates $\lambda_i$ are the same or different.

The following result can be deduced from Theorem~\ref{dcero} and Proposition \ref{otrocerogene}.
\begin{corollary}
\label{igual}
Assume that $\lambda_1=\lambda_2=\lambda$ and ${\bf \Theta}=(\theta_{\l,1},\theta_{\l,2})$ the unique solution of problem \eqref{eq:membrane.logis}. Then,
$$
\lim_{\lambda\to 0}\frac{\theta_{\l,i}}{\lambda}=\frac{|\Omega_1|+|\Omega_2|}{\displaystyle\int_{\Omega_1}\alpha_1(x)+\int_{\Omega_2}\alpha_2(x)},\quad \mbox{and}
\quad \lim_{\lambda\to \infty}\frac{\theta_{\l,i}}{\lambda}=\frac{1}{\alpha_i(x)}\quad \quad\mbox{both limits uniformly in $\overline\Omega_i.$}
$$
\end{corollary}
\begin{proof}
It is obvious that $(\theta_{\l,1}/\lambda, \theta_{\l,2}/\lambda)$ verifies (\ref{eqlogisdinfi}) with $\beta_i(x)\equiv 1$ and $d=1/\lambda$.
\end{proof}

We now focus on the more general case when the growth rates are different, i.e. $\lambda_1\neq \lambda_2$ and study the behaviour of ${\bf \Theta}_\lambda$ when $\lambda_1\to\pm\infty$ for fixed values for $\lambda_2$.

To do that, we fist introduce $L_{\lambda_2}$, the so-called large solution of the problem
\begin{equation}
\label{eq:classic_large_lambda2}
 \left\{\begin{array} {l@{\quad}l} -\Delta v =v(\lambda_2 -\alpha_2(x) v)& \hbox{in}\quad \Omega_2,\\
 v=\infty& \hbox{on}\quad \Sigma,\\
 \partial_{\bf n_2} v=0 & \hbox{on}\quad \Gamma,
\end{array}\right.
\end{equation}
 where the boundary data on $\Sigma$ has to be understood as $v(x)\to \infty$ as ${\rm dist}(x,\Sigma)\to 0$. The existence of this kind of solution is shown in~\cite{DuHuang} and~\cite{JulianBook2} and also during the proof of Proposition~\ref{Prop_v0} below.

\begin{proposition}
\label{Prop_v0}
Let $\lambda_2$ fixed such that ${\bf \Theta}_\lambda$ is the unique positive solution of problem \eqref{eq:membrane.logis}. Then,
\begin{equation}
\label{limit:lambda_1}
\lim_{\lambda_1\to \infty} \theta_{\lambda_1}(x)=+\infty,\quad \hbox{for all $x\in \overline \Omega_1$},
\end{equation}
and
\begin{equation}
\label{limit:lambda_2}
\lim_{\lambda_1\to \infty} \theta_{\l_2}=L_{\lambda_2},\quad \hbox{in $C^2(\Omega_2)$}.
\end{equation}
\end{proposition}

\begin{proof}
First, thanks to the existence results shown above, see Propositions~\ref{exislogismem} and Lemma~\ref{autovalor}, we observe that for any value of $\lambda_2$, there exists a positive solution~\eqref{eq:membrane.logis} for large $\lambda_1$.

Also, due to~(\ref{cotaimpor}), we have that
\begin{equation*}
\label{ine}
u_{\l_1}\leq \theta_{\l_1}\quad\mbox{in $\O_1$,}
\end{equation*}
where $u_{\lambda_1}$ is the solution of equation {\rm(\ref{w1-w2.sec5})}.
 Moreover, let $\varphi_1$ be a positive eigenfunction associated with $\sigma_1:=\sigma_1^{\Omega_1}[-\D;\mathcal{N}+\gamma_1].$ It is clear that if $\varepsilon>0$
 is chosen so that $ \varepsilon \alpha_1\|\varphi_1\|_\infty\leq \lambda_1-\sigma_1$, then $\varepsilon \varphi_1$ is a subsolution to~{(\ref{w1-w2.sec5})}. Hence
$$
\frac{\lambda_1-\sigma_1}{\alpha_1\|\varphi_1\|_\infty}\varphi_1(x)\leq u_{\lambda_1}(x)\quad\mbox{in $\Omega_1$,}
$$
and \eqref{limit:lambda_1} follows by taking the limit as $\lambda_1\to\infty$.

On the other hand, once we know that \eqref{limit:lambda_1} is achieved it is clear, by construction, that $\partial_{\bf n_2} u_2$ goes to infinity on the interface $\Sigma$. Consequently, since $u_2$ is a positive solution for its corresponding equation in
 \eqref{eq:membrane.logis} for that fixed $\lambda_2$ we actually have the limit \eqref{limit:lambda_2}. Therefore, $u_2$ is a large solution of the problem \eqref{eq:classic_large_lambda2}. We sketch the proof.

 Let  $v_m$ denote the unique positive solution of
 \begin{equation*}
\label{vm}
 \left\{\begin{array}{ll}
 -\Delta v =v(\lambda_2 -\alpha_2(x) v)& \mbox{in $\Omega_2$,}\\
 \partial_{\bf n_2} v+\g_2v=m& \hbox{on $\Sigma$,}\\
 \partial_{\bf n_2} v=0 & \hbox{on $\Gamma$.}
\end{array}\right.
\end{equation*}
 It is evident that
 \begin{equation*}
 \label{cadena}
 v_m<\theta_{\lambda_2}<v_M,
 \end{equation*}
 where $m=\gamma_2(\displaystyle\min_{x\in \Sigma}\theta_{\lambda_1}(x))$ and $M=\gamma_2(\displaystyle\max_{x\in \Sigma}\theta_{\lambda_1}(x))$.
 Furthermore, we define
 $$
 z(x):=\lim_{m\to\infty}v_m.
 $$
Since $v_m\leq L_{\l_2}$ for all $m$, we deduce that $z\leq L_{\l_2}$ in $\Omega_2$. Hence, it is enough to show that
 $$
 z(x)\to \infty\quad\mbox{as $\delta(x)\to 0$,} \quad \hbox{where}\quad \delta(x):={\rm dist}(x,\Sigma)\to 0.
 $$
To this end, we follow the ideas performed in \cite{GM-R-S2}. For $\delta_0>0$ consider the set
 $$
 \Omega_{\delta_0}:=\{x\in\Omega_2:0<\delta(x)<\delta_0\},\quad \delta_0>0,
 $$
 which is a small band close to $\Sigma$, so that we take $\delta_0$ small enough such that $\delta \in C^2(\overline\Omega_\delta)$ and $|\nabla \d(x)|=1$ in $\overline\Omega_{\delta_0}$.

Also, for $A>0$, $k>0$ and $\tau>0$, we define the auxiliary function
 $$
 \underline{w}(x):=A(\delta(x)+\tau)^{-2}-k.
 $$
and for $\Gamma_{\delta_0}:=\{x\in\Omega_2\,:\,   \delta(x)=\delta_0\}$ we  consider the problem
\begin{equation}
\label{largam}
\left\{
\begin{array}{ll}
 -\Delta v=v(\lambda_2-\alpha_2(x)v) \quad & \mbox{in $\Omega_{\delta_0}$,}\\
 \partial_{\bf n_2} v+\gamma_2v=m &  \mbox{on $\Sigma$,}\\
 v=v_m &  \mbox{on $\Gamma_{\delta_0}$.}
\end{array}
\right.
\end{equation}
 Thus, $\underline{w}$ is subsolution of (\ref{largam}) if the following three conditions hold:
 \begin{equation}
 \label{c1}
 \begin{aligned}
 2A(\delta(x)&+\tau)^{-4}\left\{-3|\nabla \delta|^2+\alpha_2\frac{A}{2}+(d(x)+\tau)\Delta \delta- (\delta(x)+\tau)^{2}\frac{2k\alpha_2+\lambda_2}{2}\right\}\\
&+k(\alpha_2 k+\l_2)\leq 0\quad\mbox{in $\Omega_{\delta_0}$,}
 \end{aligned}
 \end{equation}
 \begin{equation}
 \label{c2}
  -2A\tau^{-3}\partial_{\bf n_2}  \delta+\gamma_2(A \tau^{-2}-k)\leq m \quad\mbox{on $\Sigma$,}
 \end{equation}
 and
 \begin{equation}
 \label{c3}
 A({\delta_0}+\tau)^{-2}-k\leq v_m\quad\mbox{on $\Gamma_{\delta_0}$.}
  \end{equation}
 Analyzing  (\ref{c3}), since $v_m$ is increasing in $m$, (\ref{c3}) holds if
 \begin{equation}
 \label{eq:expresion_k}
 k= k(\delta_0, \tau):=A(\delta_0+\tau)^{-2}-V_1(\delta_0),
 \end{equation}
 where $V_1(\delta):=\max_{\Gamma_{\delta_0}}v_1(x)$.
 Now, let us study (\ref{c1}).  Observe that
$$
 (\delta(x)+\tau)\Delta \delta
- (\delta(x)+\tau)^{2}\frac{2k\alpha_2+\l_2}{2}  <C(\delta(x)+\tau)- (\delta(x)+\tau)^{2}\frac{2k\alpha_2+\lambda_2}{2} <\frac{C^2}{2(2k\alpha_2+\lambda_2)},
 $$
 where $C=\max|\Delta \delta|$.
 Moreover, since $|\nabla \delta|=1$ we get to
 $$
 -3|\nabla \delta|^2+\alpha_2\frac{A}{2}+(\delta(x)+\tau)\Delta \delta
- (\delta(x)+\tau)^{2}\left[\frac{2k\alpha_2+\lambda_2}{2}\right]<-3+\alpha_2\frac{A}{2}+\frac{C^2}{2(2k\alpha_2+\lambda_2)}.
$$
 Taking the expression \eqref{eq:expresion_k} of $k$, we get that (\ref{c1}) holds if
 $$
 2A\left(-3+\alpha_2\frac{A}{2}+\frac{C^2}{2(2k\alpha_2+\lambda_2)}\right)+\alpha_2A^2<0,
 $$
 or equivalently,
 $$
 -6+2\alpha_2A+\frac{C^2}{2k\alpha_2+\lambda_2}<0.
 $$
 Finally, choosing $A<3/\alpha_2$, we can take $\delta_0$ and $\tau$ small enough such that the above inequality is verified.
 Moreover, for the value of $k$ given by \eqref{eq:expresion_k}, we take $m$ large such that (\ref{c2}) holds.

Consequently for $\delta_0>0$ sufficiently small, there exist $\tau_0=\tau_0(\delta_0)$ and $k=k(\delta_0)$ such that $\underline w$ is subsolution of problem (\ref{largam}) for $m\geq m(\tau)$. Hence,
 $$
 A(\delta(x)+\tau)^{-2}-k\leq v_m(x)\quad x\in \Omega_{\delta_0}.
 $$
 Therefore, passing to the limit as $m\to\infty$ and $\tau\to 0$, we get
 $$
  A\delta(x)^{-2}-k\leq z(x),
 $$
 whence we deduce that $z\to\infty$ as $\delta(x)\to 0$.
\end{proof}

Finally we analyze the behaviour of the solution ${\bf \Theta}_\lambda$ of~\eqref{eq:membrane.logis} as $\l_1\to -\infty$.

\begin{proposition}
\label{minusinfu}
If  $\lambda_2\geq \sigma_2$, then,  as $\lambda_1\to-\infty$,
\begin{equation*}
\label{muynega}
\theta_{\lambda_1}\to 0\quad\mbox{in $L^\infty(\Omega_1)$ \quad and }\quad
\theta_{\lambda_2}\to w_2\quad\mbox{in $C^2(\overline\Omega_2)$. }
\end{equation*}
On the other hand,  if $\lambda_2<\sigma_2$, then, there is no positive solution of~{\rm \eqref{eq:membrane.logis}} for a very negative $\lambda_1$.
\end{proposition}

\begin{proof}
Note first that the second statement ($\lambda_2<\sigma_2$) follows directly from Remark  \ref{dis}. Analogously, if $\lambda_2\geq \sigma_2$  from Remark  \ref{dis} follows that ${\bf \Theta_\lambda}$ exists and its unique.

To analyze the behaviour of ${\bf \Theta_\lambda}$ when $\lambda_2\geq \sigma_2$, observe that $\theta_{\lambda_1}\leq z_1$ where $z_1$ is the unique positive solution of
 \begin{equation}
\label{z1}
\left\{
\begin{array}{ll}
 -\Delta z_1=\lambda_1 z_1 & \mbox{in $\Omega_1$,}\\
  \partial_{\bf n_1} z_1+\gamma_1z_1=\gamma_1 K &  \mbox{on $\Sigma$,}
\end{array}
\right.
\end{equation}
where $K=\max\{ \frac{(\beta_1)_M}{(\alpha_1)_L}, \frac{(\beta_2)_M}{(\alpha_2)_L}\} $, which exists and it is unique because $\sigma_1^{\Omega_1}[-\Delta-\lambda_1;\mathcal{N}+\gamma_1]>0$ for a very negative $\lambda_1$ (see \cite{Santi2}).

Now, we claim that
\begin{equation}
\label{am}
z_1(x)\leq C(-\lambda_1)^{-1/2}\quad \mbox{for $-\lambda_1$ large.}
\end{equation}
To prove the claim we use the blow-up argument described in \cite{GS}. Assume, by contradiction, that (\ref{am}) does not hold. Then there exist two sequences
$(-\lambda_{1,n})\nearrow\infty$, and $z_n\in C^{2,\a}(\overline\Omega_1)$, where $z_n$ represent the solutions of problem
(\ref{z1}) with $\lambda_1=\lambda_{1,n}$ such that
\begin{equation*}
\label{contra}
(-\lambda_{1,n})^{1/2} M_n\nearrow \infty,  \quad \hbox{where}\quad M_n=\|z_n\|_\infty.
\end{equation*}
Let  $x_n \in
\overline\Omega_1$ be a point where $z_n$ attains its maximum and assume with no loss of
generality that $x_n\to x_0\in \overline\Omega_1$. Thus, to complete the proof we need to distinguish two different
cases: either $x_0\in \Omega_1$ or $x_0\in \Sigma$.

\smallskip

\noindent {\it Case 1.} If $x_0\in \Omega_1$ we introduce the scaled
functions
$$
v_n (y)= \frac{z_n(x_n+(-\lambda_{1,n})^{-1/2} y)}{M_n},
$$
which verify that $v_n(0)=1$, $0\le v_n\le 1$  and
$$
-\Delta v_n+v_n=0  \qquad \hbox{in }
\Omega_n:=(-\lambda_{1,n})^{1/2} (-x_n +\Omega_1).
$$
Moreover, given $R>0$, there exists $n_0\in\mathbb{N}$ such that $x_n+(-\lambda_{1,n})^{-1/2} y\in\Omega_1$ for $n\geq n_0$ and $y\in B_R(0)$.
Hence, $v_n$ is bounded in $W^{2,p}(B_R(0))$ for all $p>1$. Hence $v_n\to v$ in $C^1(\overline\Omega_1)$ and $v$ verifies
$$
-\Delta v+v=0\quad\mbox{in $B_R(0)$ for each $R>0$,}
$$
with $0\le v \le 1$, $v(0)=1$, which is impossible by the maximum principle.
\smallskip

\noindent {\it Case 2.} If $x_0\in \Sigma$, under a new change of variable we can assume that a neighbourhood
around $x_0\in \Sigma$ is contained in the hyperplane $\{x^N=0\}$ and, also, that the domain $\overline\Omega_1$ is contained  in the set $H:=\{x\in \mathbb{R}^N:x^N>0\}$. Then, given $R>0$ for $n\geq n_0$, $v_n$ is defined in
$$
H_{R,n}:=B_R(0)\cap \left\{y^N>-\frac{x_n^N}{(-\lambda_1)_n^{1/2}}\right\}.
$$
Observe that since $x_n^N$ is bounded,
we have ${x_n^N}{(-\l_{1,n})^{-1/2}}\to 0$, as $n\to \infty$. Thus, the set $H_{R,n}$ tends to
 $B_R(0) \cap H$, as $n\to \infty$. In this particular case, $v_n$ verifies
\begin{equation*}
\label{case2.1}
\left\{\begin{array}{ll}
-\Delta v_n+v_n =0   & \mbox{in $H\cap B_R$,}\\
\displaystyle \partial_{\bf n_1}v_n+\gamma_1 \frac{v_n}{(-\lambda_{1,n})^{1/2}} =\frac{\gamma_1K }{M_n(-\lambda_{1,n})^{1/2}} & \mbox{on $\partial H\cap B_R$.}
\end{array}\right.
\end{equation*}
Therefore, passing to the limit we arrive at the existence of a function $v\in C^2(\overline H)$,
which is a solution of the problem
\begin{equation*}\label{eq4}
\left\{\begin{array}{ll}
-\Delta v+v=0 & \hbox{in  $H$,}\\
\displaystyle\frac{\partial v}{\partial {\bf n}}=0 & \hbox{on $\partial H$,}
\end{array}\right.
\end{equation*}
with $0\le v \le 1$, $v(0)=1$. Applying the reflection principle, we have that $v\in C^2(\overline H)$, with $v(0) = 1$ and such that $v$
is a solution of $-\Delta v +v=0$ in $\mathbb{R}^N$, a contradiction.
\end{proof}


\end{document}